\documentstyle[amscd,amssymb,verbatim]{amsart}
\newcommand{\lrar}[1]{\begin{picture}(50,10)(-25,-5)                          
\put(-25,0){\vector(1,0){50}}
\put(0,5){\makebox(0,0)[b]{\mbox{$#1$}}}
\end{picture}}

\newcommand{\ev}{\operatorname{ev}}
\newcommand{\ad}{\operatorname{ad}}

\newcommand{\pr}{\operatorname{pr}}

\newcommand{\splin}{\operatorname{sl}}

\newcommand{\hg}{{\frak h}}

\newcommand{\lan}{\langle}
\newcommand{\ran}{\rangle}

\newcommand{\Mat}{\operatorname{Mat}}

\newcommand{\tr}{\operatorname{tr}}

\renewcommand{\P}{{\Bbb P}}

\newcommand{\si}{\sigma}

\newcommand{\Pic}{\operatorname{Pic}}

\newcommand{\ga}{\gamma}
\newcommand{\Si}{\Sigma}
\newcommand{\de}{\delta}
\newcommand{\eps}{\epsilon}
\newcommand{\We}{\bigwedge}

\numberwithin{equation}{section}

\newtheorem{thm}{Theorem}[section]
\newtheorem{prop}[thm]{Proposition}
\newtheorem{lem}[thm]{Lemma}

\newenvironment{defi}{\vspace{3mm}\noindent
{\bf Definition.}}{\vspace{3mm}}
\newenvironment{rem}{\vspace{3mm}\noindent
{\bf Remark.}}{\vspace{3mm}}
\newenvironment{rems}{\vspace{3mm}
\noindent {\bf Remarks.}}{\vspace{3mm}}
\newenvironment{ex}{\vspace{3mm}\noindent
{\bf Example.}}{\vspace{3mm}}

\newcommand{\Pf}{\noindent {\it Proof}}
\newcommand{\id}{\operatorname{id}}

\newcommand{\ov}{\overline}
\newcommand{\pa}{\partial}

\newcommand{\ra}{\rightarrow}

\newcommand{\OO}{{\cal O}}

\newcommand{\Hom}{\operatorname{Hom}}
\newcommand{\Ext}{\operatorname{Ext}}
\newcommand{\End}{\operatorname{End}}
\newcommand{\Res}{\operatorname{Res}}

\renewcommand{\a}{\alpha}
\renewcommand{\b}{\beta}
\newcommand{\bfe}{{\bf e}}
\newcommand{\om}{\omega}

\newcommand{\la}{\lambda}
\newcommand{\th}{\theta}
\newcommand{\C}{{\Bbb C}}

\newcommand{\Z}{{\Bbb Z}}

\newcommand{\bm}{{\bf m}}
\newcommand{\La}{\Lambda}
\newcommand{\Ga}{\Gamma}

\newcommand{\wt}{\widetilde}
\newcommand{\ot}{\otimes}

\newcommand{\sub}{\subset}
\newcommand{\ed}{\qed\vspace{3mm}}
\newcommand{\const}{\operatorname{const}}
\newcommand{\GL}{\operatorname{GL}}

\title{Massey products on cycles of projective lines and trigonometric solutions of the 
Yang-Baxter equations}
\author{A. Polishchuk}
\address{Department of Mathematics, University of Oregon, Eugene, OR 97403}
\email{apolish@@uoregon.edu}
\thanks{This work was partially supported by the NSF grant DMS-0601034}

\begin{document}

\begin{abstract} We show that a nondegenerate unitary solution $r(u,v)$ of
the associative Yang-Baxter equation (AYBE) for $\Mat(N,\C)$ (see \cite{P-AYBE}) with the Laurent series at $u=0$ of the form $r(u,v)=\frac{1\ot 1}{u}+r_0(v)+\ldots$
satisfies the quantum Yang-Baxter equation, provided
the projection of $r_0(v)$ to $\splin_N\ot\splin_N$ has a period.
We classify all such solutions of the AYBE extending the work of Schedler \cite{Sch}. We also characterize solutions coming from
triple Massey products in the derived category of coherent sheaves on cycles of projective lines.
\end{abstract}

\maketitle

\bigskip

\centerline{\sc Introduction}

\bigskip

This paper is concerned with solutions of the associative Yang-Baxter equation (AYBE)
\begin{equation}\label{AYBE}
r^{12}(-u',v)r^{13}(u+u',v+v')-
r^{23}(u+u',v')r^{12}(u,v)+
r^{13}(u,v+v')r^{23}(u',v')=0,
\end{equation}
where $r(u,v)$ is a meromorphic function of two complex variables $(u,v)$
in a neighborhood of $(0,0)$ taking values in $A\otimes A$, where
$A=\Mat(N,\C)$ is the matrix algebra. Here we use the notation $r^{12}=r\ot 1\in A\ot A\ot A$, etc.
We will refer to a solution of \eqref{AYBE} as an {\it associative $r$-matrix}.
This equation was introduced in the above form in 
\cite{P-AYBE} in connection with triple Massey products for simple vector bundles on elliptic curves and their degenerations. It is usually coupled with the {\it unitarity} condition 
\begin{equation}\label{skewsymcond}
r^{21}(-u,-v)=-r(u,v).
\end{equation}
Note that 
the constant version of \eqref{AYBE} was independently introduced in \cite{A} in connection with the
notion of infinitesimal bialgebra (where $A$ can be any associative algebra).
The AYBE is closely related to the classical Yang-Baxter equation (CYBE) with spectral parameter
\begin{equation}\label{CYBE}
[r^{12}(v),r^{13}(v+v')]-[r^{23}(v'),r^{12}(v)]+[r^{13}(v+v'),r^{23}(v')]=0
\end{equation}
for the Lie algebra $\splin_N$ (so $r(v)$ takes values in $\splin_N\ot\splin_N$) and also with the quantum Yang-Baxter equation (QYBE)
with spectral parameter
\begin{equation}\label{QYBE}
R^{12}(v)R^{13}(v+v')R^{23}(v')=R^{23}(v')R^{13}(v+v')R^{12}(v),
\end{equation}
where $R(v)$ takes values in $A\ot A$.
In the seminal work \cite{BD} Belavin and Drinfeld made a thorough study of
the CYBE for simple Lie algebras. In particular, they showed that all nondegenerate solutions are equivalent to either elliptic, trigonometric, or rational solutions, and gave a complete classification in the elliptic and trigonometric cases. In the present paper we extend some of their results and techniques to the AYBE. In addition, we show that often solutions of the AYBE are automatically solutions
of the QYBE (for fixed $u$).

We will be mostly studying unitary solutions of the AYBE (i.e., solutions of \eqref{AYBE} and
\eqref{skewsymcond}) that have the Laurent expansion at $u=0$ of the form
\begin{equation}\label{Laur-exp}
r(u,v)=\frac{1\ot 1}{u}+r_0(v)+ur_1(v)+\ldots
\end{equation}
It is easy to see that in this case $r_0(v)$ is a solution of the CYBE. Hence, denoting
by $\pr:\Mat(N,\C)\to\splin_N$ the projection along $\C\cdot 1$ we obtain that 
$\ov{r}_0(v)=(\pr\ot\pr)r_0(v)$ is a solution of the CYBE for $\splin_N$. 
We prove that if $r(u,v)$ is nondegenerate (i.e., 
the tensor $r(u,v)\in A\ot A$ is nondegenerate for generic $(u,v)$) then
so is $\ov{r}_0$. Thus, $\ov{r}_0$ falls within Belavin-Drinfeld classification.
Furthermore, we show that if $\ov{r}_0$ is either elliptic or trigonometric then $r(u,v)$ is uniquely determined by $\ov{r}_0$ up to certain natural transformations. 
The natural question raised in \cite{P-AYBE} is which solutions of the CYBE for $\splin_N$
extend to unitary solutions of the AYBE of the form \eqref{Laur-exp}. 
In \cite{P-AYBE} we showed that this is the case for all elliptic solutions and gave some examples
with trigonometric solutions. In \cite{Sch} Schedler studied further this question for trigonometric
solutions of the CYBE of the form $r_0(v)=\frac{r+e^v r^{21}}{1-e^v}$, where $r$ is a constant solution
of the CYBE. He discovered that not all trigonometric solutions of the CYBE can be extended to
solutions of the AYBE, and found a nice combinatorial structure that governs the situation 
(called {\it associative BD triples}). In this paper we complete the picture by giving the answer
to the above question for arbitrary trigonometric solutions of the CYBE (see Theorem \ref{main-thm}
below). We will also prove that every 
nondegenerate unitary solution $r(u,v)$ of the AYBE with the Laurent expansion at $u=0$ of the form
\eqref{Laur-exp} satisfies the QYBE with spectral parameter
for fixed $u$, provided $\ov{r}_0(v)$ either has a period (i.e., it is either elliptic or trigonometric) or has 
no infinitesimal symmetries (see Theorem \ref{QYBE-thm}). Thus, our work on extending
trigonometric classical $r$-matrices (with spectral parameter) to solutions of the AYBE leads to explicit formulas for the corresponding quantum $r$-matrices. The connection with the QYBE was noticed before for elliptic solutions constructed in \cite{P-AYBE} 
(because they are given essentially by Belavin's elliptic $R$-matrix) and 
also for those trigonometric solutions that are constructed in \cite{Sch}. 

An important input for our study of trigonometric solutions of the AYBE is the geometric picture
with Massey products developed in \cite{P-AYBE}
that involves considering simple vector bundles on elliptic curves and their
rational degenerations. In {\it loc. cit.} we constructed all elliptic
solutions in this way and some trigonometric solutions coming from simple vector bundles on the union of two projective lines glued at two points. In this paper we consider the case of bundles on a cycle of projective lines of arbitrary length. We compute explicitly corresponding solutions of the AYBE. Then we notice that similar formula make sense in a more general context and prove this by a direct calculation.
The completeness of the obtained list of trigonometric solutions is then checked by combining the
arguments of \cite{Sch} with those of \cite{BD} (modified appropriately for the case of the AYBE).
It is interesting that contrary to the initial expectation expressed in \cite{P-AYBE} not all trigonometric solutions of the AYBE can be obtained from the triple Massey products on cycles of projective lines (see Theorem \ref{geom-BD-thm}). This makes us wonder whether there is some generalization
of our geometric setup. 

Another question that seems to be worth pursuing is the connection between the combinatorics
of simple vector bundles on a cycle of projective lines $X$ and the Belavin-Drinfeld combinatorics.
Namely, the discrete type of a vector bundle on $X$ 
is described by the splitting type on each component of $X$. As was observed in 
\cite{BDG}, Theorem 5.3, 
simplicity of a vector bundle corresponds to a certain combinatorial condition on these splitting types
(see also Lemma \ref{simple-crit-lem}). In this paper we show that this condition allows to associate
with such a splitting type a Belavin-Drinfeld triple (or rather an enhanced combinatorial data described
below). It seems that this connection might provide an additional insight on the problem of classifying discrete types of simple vector bundles on $X$.

In \cite{Mud} Mudrov constructs solutions of the QYBE from certain algebraic data that should be viewed as associative analogues of Manin triples. Elsewhere we will show how solutions of the AYBE
give rise to such data and will study the corresponding associative algebras that are related to
both the classical and quantum side of the story.

Now let us present the combinatorial data on which our trigonometric solutions of the AYBE depend
(generalizing Belavin-Drinfeld triples with associative structure considered in \cite{Sch}).
Let $S$ be a finite set. To equip $S$ with a cyclic order is the same as to fix a transitive cyclic permutation $C_0:S\to S$. We denote by $\Ga_{C_0}:=\{(s,C_0(s))\ |\ s\in S\}$ the graph of $C_0$. 

\begin{defi} An {\it associative BD-structure} on a finite set 
$S$ is given by a pair of transitive cyclic permutations $C_0,C:S\to S$ and
a pair of proper subsets $\Ga_1,\Ga_2\sub\Ga_{C_0}$, such that
$(C\times C)(\Ga_1)=\Ga_2$, where $(C\times C)(i,i')=(C(i),C(i'))$.
\end{defi}

We can identify $\Ga_{C_0}$ with the set of vertices $\Ga$ of the affine Dynkin diagram 
$\wt{A}_{N-1}$, where $N=|S|$ (preserving the cyclic order). Then
we get from the above structure a Belavin-Drinfeld triple $(\Ga_1,\Ga_2,\tau)$
for $\wt{A}_{N-1}$, where the bijection $\tau:\Ga_1\to\Ga_2$ is induced by $C\times C$.
It is clear that $\tau$ preserves the inner product. The nilpotency condition on $\tau$
is satisfied automatically. Indeed, choose $(s_1,C_0(s_1))\in\Ga_S\setminus\Ga_1$. Then
for every $(s,C_0(s))\in\Ga_1$ there exists $k\ge 1$ with $C^k(s)=s_1$, so that
$(C\times C)^k(s,C_0(s))\not\in\Ga_1$.

We extend the bijection $\tau$ to a bijection $\tau:P_1\to P_2$ induced by $C\times C$,
where 
$$P_{\iota}=\{(s,C_0^k(s))\ |\ (s,C_0(s))\in\Ga_{\iota}, (C_0(s),C_0^2(s))\in\Ga_{\iota},\ldots, 
(C_0^{k-1}(s),C_0^k(s))\in\Ga_{\iota}\},\ \iota=1,2.$$

For a finite set $S$ let us denote by $A_S$
the algebra of endomorphisms of the $\C$-vector space with the basis $(\bfe_i)_{i\in S}$,
so that $A_S\simeq\Mat(N,\C)$,  where $N=|S|$. We denote by $e_{ij}\in A_S$ the endomorphism
defined by $e_{ij}(\bfe_k)=\de_{jk}\bfe_i$.
We denote by $\hg\sub A_S$ the subalgebra of diagonal matrices (i.e., the span of $(e_{ii})_{i\in S}$).
Now we can formulate our result about trigonometric solutions of the AYBE.

\begin{thm}\label{main-thm}
(i) Let $(C_0,C,\Ga_1,\Ga_2)$ be an associative BD-structure on a finite set $S$.
Consider the $A_S\ot A_S$-valued function 
\begin{align*}
&r(u,v)=\frac{1}{1-\exp(-v)}\sum_i e_{ii}\ot e_{ii}+\frac{1}{\exp(u)-1}\sum_{0\le k<N, i}
\exp(\frac{ku}{N})e_{C^k(i),C^k(i)}\ot e_{ii}+\\
&\frac{1}{\exp(v)-1}\sum_{0<m<N, j=C_0^m(i)}\exp(\frac{mv}{N})e_{ij}\ot e_{ji}+\\
&\sum_{0<m<N,k\ge 1;j=C_0^m(i),\tau^k(i,j)=(i',j')}
[\exp(-\frac{ku+mv}{N})e_{ji}\ot e_{i'j'}-\exp(\frac{ku+mv}{N})e_{i'j'}\ot e_{ji}],
\end{align*}
where $i,i',j,j'$ denote elements of $S$, and the summation in the last sum is taken only over
those $(i,j)$ for which $\tau^k$ is defined on $(i,j)$. 
Then $r(u,v)$ satisfies \eqref{AYBE} and \eqref{skewsymcond}. 
Furthermore, let us set
\begin{equation}\label{Q-r-matrix}
R(u,v)=\left([\exp(\frac{u}{2})-\exp(-\frac{u}{2})]^{-1}+[\exp(\frac{v}{2})-\exp(-\frac{v}{2})]^{-1}\right)^{-1}\cdot
r(u,v).
\end{equation}
Then $R(u,v)$ satisfies the QYBE with spectral parameter \eqref{QYBE} (for fixed $u$) and the unitarity
condition
\begin{equation}\label{QYBE-un-eq}
R(u,v)R^{21}(u,-v)=1\ot 1.
\end{equation}

\noindent
(ii) Assume that $N>1$. Then every nondegenerate unitary solution of the AYBE for $A=\Mat(N,\C)$
with the Laurent expansion at $u=0$ of
the form \eqref{Laur-exp}, where $\ov{r}_0(v)$ is a trigonometric solution of the CYBE
for $\splin_N$, is equal to 
$$c\exp(\la uv)\exp[u(1\ot a)+v(b\ot 1)]r(cu,c'v)\exp[-u(a\ot 1)-v(b\ot 1)],$$
where $r(u,v)$ is obtained from one of the solutions from (i) by applying an algebra isomorphism
$A_S\simeq A$, $\la$, $c$ and $c'$ are constants ($c\neq 0$, $c'\neq 0$),
and $a,b\in\hg$ are infinitesimal symmetries of $r(u,v)$, i.e., 
$$[a\ot 1+1\ot a,r(u,v)]=[b\ot 1+1\ot b,r(u,v)]=0.$$
\end{thm}

Note that the complete list of scalar unitary solutions of the AYBE was obtained in Theorem 5 of
\cite{P-AYBE}. The solution obtained from Theorem \ref{main-thm}(i) in the case $N=1$ 
coincides with the basic trigonometric solution from that list (up to changing $v$ to $-v$).

We will also deduce the following result about solutions of the AYBE not depending on the
variable $u$.

\begin{thm}\label{AYBE-v} Assume that $N>1$. 
Let $r(v)$ be a nondegenerate unitary solution of the AYBE for $A=\Mat(N,\C)$ not depending
on the variable $u$. Then 
$$r(v)=\ov{r}(v)+b\ot 1+1\ot b+\frac{c\cdot 1\ot 1}{Nv},$$
where $\ov{r}(v)$ is equivalent to a rational nondegenerate solution of the CYBE for $\splin_N$, 
$b\in\splin_N$ is an infinitesimal symmetry of $\ov{r}(v)$, $c\in\C^*$. 
Also, 
$$R(u,v)=\left(1+\frac{cu}{v}\right)^{-1}\cdot\left(1+ur(v)\right)$$
is a unitary solution of the QYBE with spectral parameter for fixed $u$
(hence, the same is true for $\frac{v}{c}r(v)=\lim_{u\to\infty}R(u,v)$).
\end{thm} 

The case of nondegenerate unitary solutions of the AYBE not depending on $v$ turns out to be much easier --- in this case we get a complete list of solutions (see Proposition \ref{u-prop}). Note that
there are no constant nondegenerate solutions of the AYBE for $A=\Mat(N,\C)$ (unitary or not),
as follows from Proposition 2.9 of \cite{A2}. 

The paper is organized as follows.
In section \ref{QYBE-sec} we discuss nondegeneracy
conditions for solutions of the AYBE and show how to deduce the QYBE in Theorem \ref{QYBE-thm}.
After recalling in section \ref{bundle-sec}
the geometric setup leading to solutions of the AYBE, we calculate these solutions associated
with simple vector bundles on cycles of projective lines in sections \ref{cycles-sec} and
\ref{comp-sec} (the result is given by formulas \eqref{r-const-eq}, \eqref{r-gen-eq}).
Then in section \ref{comb-sec} we consider associative BD-structures on completely ordered sets and classify such structures coming from simple vector bundles on cycles of projective lines (see Theorem 
\ref{geom-BD-thm}). In section \ref{comb-sol-sec} we prove the first part of Theorem \ref{main-thm}.
In section \ref{pole-sec} we establish a meromorphic continuation in $v$ for a class of solutions of the AYBE and derive some additional information about these solutions. Finally, in section \ref{class-sec} we prove the second part of Theorem \ref{main-thm} and Theorem
\ref{AYBE-v}.

{\it Acknowledgment}. I am grateful to Pavel Etingof for the crucial help with organizing my
initial computations into a nice combinatorial pattern. I also thank him and Travis Schedler for
useful comments on the first draft of the paper and the subsequent helpful discussions. Parts of this work were done while the author 
enjoyed the hospitality of the Max-Planck-Institute f\"ur Mathematik in Bonn and of the SISSA in Trieste.

\section{The AYBE and the QYBE}
\label{QYBE-sec}

Recall that we denote $A=\Mat(N,\C)$.
Let $r(u,v)$ be a meromorphic $A\otimes A$-function in a neighborhood of $(0,0)$.
We say that $r(u,v)$ is {\it nondegenerate} if the tensor $r(u,v)$ is nondegenerate
for generic $(u,v)$.

We start by collecting some facts about nondegenerate unitary solutions of the AYBE.
First, let us consider the case when $r(u,v)$ does not depend on $v$. Then the AYBE reduces to
\begin{equation}\label{AYBE-u}
r^{12}(-u')r^{13}(u+u')-
r^{23}(u+u')r^{12}(u)+
r^{13}(u)r^{23}(u')=0,
\end{equation}
and the unitarity condition becomes $r^{21}(-u)=-r(u)$.

Let us set $P=\sum_{i,j} e_{ij}\ot e_{ji}$.

\begin{prop}\label{u-prop}
All nondegenerate unitary solutions of \eqref{AYBE-u} have form
$$r(u)=(\phi_a(cu)\otimes\id)(P),$$
where $c\in\C^*$, $a\in \splin_N$,
$\phi_a(u)\in\End(A)$ is the linear operator on $A$ defined from the equation
$$u\phi_a(u)(X)+[a,\phi_a(u)(X)]=X.$$
\end{prop}

\Pf . Let us write $r(u,v)$ in the form
$r(u)=(\phi(u)\otimes\id)(e)$, where $e\in A^*\otimes A$ is the canonical element,
$\phi(u):A^*\to A$ is an operator, nondegenerate for generic $u$.
Now set $B(u)(X,Y)=(X,\phi(u)^{-1}(Y))$ for $X,Y\in A$.
It is easy to see that the equation \eqref{AYBE-u} together with the unitarity condition are equivalent to the following equations on $B(u)$:
\begin{equation}\label{cocycle-eq}
B(-u)(XY,Z)+B(-u')(YZ,X)+B(u+u')(ZX,Y)=0,
\end{equation}
\begin{equation}\label{skew-bilin-eq}
B(u)(X,Y)+B(-u)(Y,X)=0.
\end{equation}
Substituting $Z=1$ in the first equation we find
$$B(-u)(XY,1)+(B(u+u')-B(u'))(X,Y)=0,\text{ i.e.},$$
$$B(u+u')(X,Y)=\xi(u)(XY)+B(u')(X,Y),$$
where $\xi(u)(X)=-B(-u)(X,1)$. Exchanging $u$ and $u'$ we
get that $C(X,Y)=B(u)(X,Y)-\xi(u)(XY)$ does not depend on $u$.
Substituting $B(u)(X,Y)=\xi(u)(XY)+C(X,Y)$ into the previous equation we get
$$\xi(u+u')=\xi(u)+\xi(u'),$$
hence, $\xi(u)=u\cdot\xi$ for some $\xi\in A^*$.
Now substituting $B(u)(X,Y)=u\cdot\xi(XY)+C(X,Y)$ into
\eqref{skew-bilin-eq} we derive that $\xi(XY)=\xi(YX)$ and $C$ is skew-symmetric.
Therefore, $\xi=c\cdot\tr$. Finally, equation \eqref{cocycle-eq} reduces to the equation
$$C(XY,Z)+C(YZ,X)+C(ZX,Y)=0.$$
Together with the skew-symmetry of $C$ this implies that $C(X,1)=C(1,X)=0$ and
the restriction of $C$ to $\splin_N\times\splin_N$ is a $2$-cocycle. Hence,
$C(X,Y)=l(XY-YX)$ for some linear functional $l$ on $\splin_N$. Conversely, for
$C$ of this form the above equation is satisfied.
Thus, all solutions of \eqref{cocycle-eq} and \eqref{skew-bilin-eq} are given by
$$B(u)(X,Y)=cu\tr(X,Y)+l(XY-YX),$$
where $c\in \C^*$ and $l$ is a linear functional on $\splin_N$.
Let us identify $A$ with $A^*$ using the metric $\tr(XY)$. Then we can view 
$\phi(u)$ as an operator from $A$ to $A$ such that
$B(u)(X,Y)=\tr(X\phi(u)^{-1}(Y))$.
Representing the functional $l$ in the form $l(X)=-\tr(Xa)$ we obtain the formula 
$$\phi(u)^{-1}(Y)=cuY+[a,Y].$$
\ed

\begin{rem} It is easy to see that $\phi_a(u)$ (and hence the corresponding associative $r$-matrix)
always has a pole at $u=0$ with order equal to the maximal $k$ such that there exists
$X\in A$ with $\ad^k(a)(X)=0$ and $\ad^{k-1}(a)(X)\neq 0$.
Indeed, $\phi_a(u)$ cannot be regular at $u=0$ since this would give
$[a,\phi_a(0)(1)]=1$. Let 
$$\phi_a(u)=\frac{\psi_{-k}}{u^k}+\frac{\psi_{-k+1}}{u^{k-1}}+\ldots $$
be the Laurent expansion of $\phi_a(u)$. Then we have
$$\psi_{i-1}+\ad(a)\circ\psi_i=0$$
for $i\neq 0$ and
$$\psi_{-1}+\ad(a)\circ\psi_0=\id.$$
Decomposing $\End(A)$ into generalized eigenspaces of the operator $\psi\mapsto\ad(a)\circ\psi$
we see that $\psi_{-1}$ is the component of $\id\in\End(A)$ corresponding to the zero eigenvalue.
This immediately implies our claim. For example, if $a$ is semisimple then $\phi_a(u)$ has
a simple pole at $u=0$. More precisely, taking diagonal matrix $a=\sum_i a_i e_{ii}$ we get
the associative $r$-matrix 
$$r(u)=\sum_{ij}\frac{1}{u+a_i-a_j} e_{ij}\ot e_{ji}.$$
\end{rem}

The proofs of the next two results are parallel to those of Propositions 2.2 and 2.1 in \cite{BD},
respectively.

\begin{lem}\label{nondeg-no-pole-lem} 
Let $r(u,v)$ be a nondegenerate unitary solution of the AYBE.
Assume that $r(u,v)$ does not have a pole at $v=0$. Then $r(u,0)$ is still nondegenerate,
and hence has the form described in Proposition \ref{u-prop}.
\end{lem}

\Pf . Let us fix $v_0$ such that $r(u,v)$ does not have a pole at $v=v_0$ and $r(u,v_0)$
is nondegenerate for generic $u$. Then we can
define a meromorphic function $\phi(u,v)$ with values in $\End_{\C}(A)$ by the condition
$$(\phi(u,v)\ot\id)(r(u,v_0))=r(u,v).$$
We claim that this function satisfies the identity
\begin{equation}\label{mat-end-eq} 
\phi(u+u',v)(XY)=\phi(u,v)(X)\phi(u',v)(Y),
\end{equation}
where $X,Y\in A$. Indeed, since $r(u,v)$ does not have a pole at $v=0$, substituting $v'=0$ in 
\eqref{AYBE} we get
$$r^{12}(-u',v)r^{13}(u+u',v)=r^{23}(u+u',0)r^{12}(u,v)-r^{13}(u,v)r^{23}(u',0).$$
Note that the right-hand side is obtained by applying $\phi(u,v)\otimes\id\otimes\id$ to the right-hand side for $v=v_0$. Applying the above equation for $v=v_0$ we deduce that it is equal to
$$(\phi(u,v)\otimes\id\otimes\id)(r^{12}(-u',v_0)r^{13}(u+u',v_0)).$$
On the other hand, the left-hand side can be rewritten as
$$[(\phi(-u',v)\otimes\id)r(-u',v_0)]^{12}[(\phi(u+u',v)\otimes\id)r(u+u',v)]^{13}.$$
Thus, if we write $r(u,v_0)=\sum K^{\a}(u)\ot \bfe_{\a}$, where $\bfe_{\a}$ is a basis of $A$, 
then we derive
$$\phi(-u',v)(K^{\a}(-u'))\phi(u+u',v)(K^{\b}(u+u'))=\phi(u,v)(K^{\a}(-u')K^{\b}(u+u')).$$
By nondegeneracy of $r(u,v_0)$ this implies \eqref{mat-end-eq}.
Taking $Y=1$ in this equation we obtain
\begin{equation}\label{mat-end-eq2}
\phi(u+u',v)(X)=\phi(u,v)(X)\phi(u',v)(1).
\end{equation}
Similarly, we deduce that
$$\phi(u+u',v)(Y)=\phi(u',v)(1)\phi(u,v)(Y).$$
Comparing these equation we see that $\phi(u',v)(1)$ commutes with $\phi(u,v)(X)$
for any $X\in A$. Using nondegeneracy of $r(u,v)$ we derive that $\phi(u,v)(1)=f(u,v)\cdot 1$
for some scalar meromorphic function $f(u,v)$. Furthermore, we should have
$$f(u+u',v)=f(u,v)f(u',v),$$
which implies that $f(u,v)=\exp(g(v)u)$ for some function $g(v)$ holomorphic near $v=0$.
Next, from \eqref{mat-end-eq2} we obtain that $\exp(-g(v)u)\phi(u,v)$ does not
depend on $u$. Thus, all solutions of \eqref{mat-end-eq} have form
$\phi(u,v)=\exp(g(v)u)\psi(v)$, where for every $v$
$\psi(v)$ is an {\it algebra automorphism} of $A$ or zero.
By our assumption $\phi(u,v)$ does not have a pole at $v=0$. Therefore,
$\psi(v)$ is holomorphic near $v=0$. Now we use the fact that every 
algebra endomorphism of $A$ is inner, and hence 
has determinant equal to $1$
(it is enough to check this for the conjugation with a diagonalizable matrix).
Since, $\psi(v_0)=\id$ this implies that $\det\psi(v)=1$ identically.
Therefore, $\det\psi(0)=1$ and $\phi(u,0)$ is invertible.
\ed

\begin{lem}\label{pole-lem} 
Let $r(u,v)$ be a nondegenerate unitary solution of the AYBE.
Assume that $r(u,v)$ has a pole at $v=0$. Then this pole is simple and 
$\lim_{v\to 0}vr(u,v)=cP$ for some nonzero constant $c$.
\end{lem}

\Pf . Let $r(u,v)=\frac{\th(u)}{v^k}+\frac{\eta(u)}{v^{k-1}}+\ldots$
be the Laurent expansion of $r(u,v)$ near $v=0$. 
Considering the polar parts as $v'\to 0$ (resp., $v\to 0$) in \eqref{AYBE} we get
\begin{equation}\label{res-eq1}
-\th^{23}(u+u')r^{12}(u,v)+r^{13}(u,v)\th^{23}(u')=0,
\end{equation}
\begin{equation}\label{res-eq2}
\th^{12}(-u')r^{13}(u+u',v')-r^{23}(u+u',v')\th^{12}(u)=0.
\end{equation}
Let $V\sub A$ be the minimal subspace such that $\th(u)\in V\ot A$ 
(for all $u$ where $\th(u)$ is defined). 
Then we have $r^{13}(u,v)\th^{23}(u')\in A\ot V\ot A$. Hence, from \eqref{res-eq1} we
get $\th^{23}(u+u')r^{12}(u,v)\in A\ot V\ot A$. This implies that
$r^{12}(u,v)\in A\ot A_1$, where 
$$A_1=\{a\in A\ :\ (a\otimes 1)\th(u)\in V\ot A\text{ for all }u \}.$$
By nondegeneracy we get $A_1=A$, hence $AV\sub V$. Similarly, using \eqref{res-eq2}
we derive that $VA\sub V$. Thus, $V$ is a nonzero two-sided ideal in $A$, so we have
$V=A$. Now let us prove that the order of pole $k$ cannot be greater that $1$.
Indeed, assuming that $k>1$ and considering the coefficient with $v^{1-k}$ in the expansion
of \eqref{AYBE} near $v=0$ we get
$$\eta^{12}(-u')r^{13}(u+u',v')-r^{23}(u+u',v')\eta^{12}(u)+\th^{12}(-u')\frac{\pa r^{13}}{\pa v}(u+u',v')=0.$$
Now looking at polar parts at $v'=0$ we get $\th^{12}(-u')\th^{13}(u+u')=0$ which contradicts to
the equality $V=A$ established above.
Therefore, $k=1$. Now let us look at \eqref{res-eq1} again. Let us fix $u$ and consider
the subspace
$$A(u)=\{x\in A\ :\ \th(u+u')(x\ot 1)=(1\ot x)\th(u')\text{ for all }u'\}.$$
Then from \eqref{res-eq1} we get that $r(u,v)\in A\ot A(u)$. By nondegeneracy this implies
that $A(u)=A$ for generic $u$, so we get an identity
$$\th(u+u')(x\ot 1)=(1\ot x)\th(u')$$
for all $x\in A$. Taking $x=1$ we see that $\th(u)=\th$ is constant. Finally, any tensor
$\th\in A\otimes A$ with the property $\th(x\ot 1)=(1\ot x)\th$ is proportional to $P$.
\ed

Recall that if $r(u,v)$ is a solution of the AYBE with the Laurent expansion at $u=0$ of the form
\eqref{Laur-exp} then $r_0(v)$ is a unitary solution of the CYBE (see proof of Lemma 1.2 in 
\cite{P-AYBE}, or Lemma 2.9 of \cite{Sch}).
The same is true for $\ov{r}_0(v)=(\pr\ot\pr)(r_0(v))\in\splin_N\ot\splin_N$.
We will show below that the nondegeneracy of $r(u,v)$ implies that $\ov{r}_0(v)$ is also nondegenerate,
hence it is either elliptic, trigonometric, or rational. 
The first two cases are distinguished from the third by the condition that
$\ov{r}_0(v)$ is periodic with respect to $v\mapsto v+p$ for some $p\in\C^*$.

Recall that by an infinitesimal symmetry of an $A\ot A$-valued function $f(x)$
we mean an element $a\in A$ such that $[a\ot 1+1\ot a,f(x)]=0$ for all $x$.

\begin{thm}\label{QYBE-thm} Let $r(u,v)$ be a nondegenerate unitary solution of the AYBE
with the Laurent expansion at $u=0$ of the form \eqref{Laur-exp},
and let $\ov{r}_0(v)=(\pr\ot\pr)(r_0(v))$. Then 

\noindent (i) $\ov{r}_0(v)$ is a nondegenerate unitary solution of the CYBE.

\noindent (ii) The following conditions are equivalent:

\noindent (a) $r(u,v)$ satisfies the QYBE \eqref{QYBE} in $v$ (for fixed $u$);

\noindent (b) the product $r(u,v)r(-u,v)$ is a scalar multiple of $1\ot 1$;

\noindent (c) $\frac{d}{dv}(r_0(v)-\ov{r}_0(v))$ is a scalar multiple of $1\ot 1$.

\noindent (d) $(\pr\ot\pr\ot\pr)[\ov{r}_0^{12}(v)\ov{r}_0^{13}(v+v')-\ov{r}_0^{23}(v')\ov{r}_0^{12}(v)+
\ov{r}_0^{13}(v+v')\ov{r}_0^{23}(v')]=0$.

\noindent (iii) The equivalent conditions in (ii) hold when
$\ov{r}_0(v)$ either admits a period or has no infinitesimal symmetries in $\splin_N$.  
\end{thm}

\begin{rems} 
1. In fact, our proof shows that equivalent conditions in (ii) hold under the weaker assumption that the system
$$[\ov{r}_0(v),a^1+a^2]=[\ov{r_0}(v),b^1+b^2+va^1]=[b,a]=0$$
on $a,b\in\splin_N$ implies that $a=0$. 

\noindent 2. Note that the implication (b)$\implies$(a) in part  (ii) of the theorem
holds for any unitary solution of the AYBE (as follows easily from Lemma \ref{main-AYBE-lem1} below).
It is plausible that one can check condition (b) in other situations than those considered in the above
theorem. For example, we have nondegenerate unitary solutions of the AYBE of the form
$$r(u,v)=\frac{\om}{u^n}+\frac{P}{v},$$
where $n\ge 1$, $\om\in A\ot A$ satisfies $\om^{12}\om^{13}=0$ and $\om^{21}=(-1)^{n-1}\om$.
It is easy to see that these solutions satisfy $r(u,v)r(-u,v)=1\ot 1/v^2$, so they also satisfy the QYBE.
On the other hand, the solutions of the AYBE constructed in Proposition \ref{u-prop}
do not satisfy the QYBE in general.
\end{rems}

\begin{lem}\label{pole-lem2} 
Assume that $N>1$. 
Let $r(u,v)$ be a nondegenerate unitary solution of the AYBE 
with the Laurent expansion at $u=0$ of the form \eqref{Laur-exp}. Then
$r(u,v)$ has a simple pole at $v=0$ with
the residue $c\cdot P$, where $c\in\C^*$.
\end{lem}

\Pf . By Lemma \ref{pole-lem} we only have to rule out the possibility that $r(u,v)$ has no pole at $v=0$.
Assume this is the case. Then $r(u,0)$ is the solution of \eqref{AYBE-u} that has a simple pole at $u=0$
with the residue $1\ot 1$. Let $\phi(u):A\to A$ be the linear operator such that
$r(u,0)=(\phi(u)\ot\id)(P)$. Then $\phi(u)$ has the Laurent expansion at $u=0$ of the form
$$\phi(u)(X)=\frac{\tr(X)\cdot 1}{u}+\psi(X)+\ldots$$
for some operator $\psi:A\to A$. By Lemma \ref{nondeg-no-pole-lem} and Proposition \ref{u-prop}, we have
$$cu\phi(u)(X)+[a,\phi(u)(X)]=X$$
for some $c\in\C^*$ and $a\in\splin_N$. Considering the constant terms of the expansions at $u=0$
we get
$$c\tr(X)\cdot 1+[a,\psi(X)]=X.$$
It follows that $[a,\psi(X)]=X$ for all $X\in\splin_N$. 
Hence, the operator $\pr\psi|_{\splin_N}:\splin_N\to\splin_N$ is invertible. Taking in the above equality
$X\in\splin_N$ such that $\pr\psi(X)=a$ we derive that $a=0$ which leads to a contradiction.
\ed

The next two lemmas constitute the core of the proof of Theorem \ref{QYBE-thm}.

\begin{lem}\label{main-AYBE-lem1} 
For a triple of variables $u_1,u_2,u_3$ (resp., $v_1,v_2,v_3$) set
$u_{ij}=u_i-u_j$ (resp., $v_{ij}=v_i-v_j$).
Then for every unitary solution of the AYBE one has
\begin{align*}
&r^{12}(u_{12},v_{12})r^{13}(u_{23},v_{13})r^{23}(u_{12},v_{23})-
r^{23}(u_{23},v_{23})r^{13}(u_{12},v_{13})r^{12}(u_{23},v_{12})=\\
&s^{23}(u_{23},v_{23})r^{13}(u_{13},v_{13})-r^{13}(u_{13},v_{13})s^{23}(u_{21},v_{23})=\\
&r^{13}(u_{13},v_{13})s^{12}(u_{32},v_{12})-s^{12}(u_{12},v_{12})r^{13}(u_{13},v_{13}),
\end{align*}
where $s(u,v)=r(u,v)r(-u,v).$
\end{lem}

\Pf . In the following proof we will use the short-hand notation $r^{ij}(u)$ for $r^{ij}(u,v_{ij})$.
The AYBE can be rewritten as
\begin{equation}\label{AYBE-c1}
r^{12}(u_{12})r^{13}(u_{23})-r^{23}(u_{23})r^{12}(u_{13})+r^{13}(u_{13})r^{23}(u_{21})=0.
\end{equation}
On the other hand, switching indices $1$ and $2$ and using the unitarity condition we obtain
\begin{equation}\label{AYBE-c2}
r^{23}(u_{23})r^{13}(u_{12})-r^{12}(u_{12})r^{23}(u_{13})+r^{13}(u_{13})r^{12}(u_{32})=0.
\end{equation}
Multiplying \eqref{AYBE-c2} with $r^{12}(u_{23})$ on the right we get
$$
r^{23}(u_{23})r^{13}(u_{12})r^{12}(u_{23})
-r^{12}(u_{12})r^{23}(u_{13})r^{12}(u_{23})+r^{13}(u_{13})s^{12}(u_{32})=0.
$$
On the other hand, switching $u_1$ and $u_2$ 
in \eqref{AYBE-c1} and multiplying the obtained equation with $r^{12}(u_{12})$ on the left we obtain
$$s^{12}(u_{12})r^{13}(u_{13})-r^{12}(u_{12})r^{23}(u_{13})r^{12}(u_{23})+
r^{12}(u_{12})r^{13}(u_{23})r^{23}(u_{12})=0.$$
Taking the difference between these cubic equations gives
$$
r^{23}(u_{23})r^{13}(u_{12})r^{12}(u_{23})-
r^{12}(u_{12})r^{13}(u_{23})r^{23}(u_{12})=s^{12}(u_{12})r^{13}(u_{13})-r^{13}(u_{13})s^{12}(u_{32}).
$$
The other half of the required equation is obtained by switching the indices $1$ and $3$ and using
the unitarity condition.
\ed


\begin{lem}\label{main-AYBE-lem2} 
Let $r(u,v)$ be a unitary solution of the AYBE with the Laurent expansion \eqref{Laur-exp} at $u=0$.
Assume also that $r(u,v)$ has a simple pole at $v=0$ with the residue $cP$. Then one has
$$s(u,v)=r(u,v)r(-u,v)=a\ot 1+1\ot a+ (f(u)+g(v))\cdot 1\ot 1$$
with
$$g(v)=-\frac{c}{N}(\tr\ot\tr)(\frac{dr_0(v)}{dv}),$$
$$f(u)=\frac{1}{N}\tr\mu(\frac{\pa r(u,0)}{\pa u}),$$
$$a=\pr\mu(\frac{\pa r(u,0)}{\pa u}),$$
where $\mu:A\ot A\to A$ denotes the product.
Furthermore,  $a\in\splin_N$ is an infinitesimal symmetry of $r(u,v)$, and
if we write
$$r_0(v)=\ov{r}_0(v)+\a(v)\ot 1-1\ot\a(-v)+h(v)\cdot 1\ot 1,$$
where $\ov{r}_0(v)\in\splin_N\ot\splin_N$ and $\a(v)\in\splin_N$, then
$$\a(v)=\a(0)-\frac{v}{cN}a.$$
\end{lem}

\Pf . Let us write $r(u,v)=\frac{cP}{v}+\wt{r}(u,v)$, where $\wt{r}(u,v)$ does not have a pole at $v=0$.
Then we can rewrite the AYBE as follows (where $v_{ij}=v_i-v_j$):
\begin{align*}
&r^{13}(u,v_{13})r^{23}(-u+h,v_{23})=
r^{23}(h,v_{23})r^{12}(u,v_{12})-r^{12}(u-h,v_{12})r^{13}(h,v_{13})=\\
&r^{23}(h,v_{23})r^{12}(u,v_{12})-r^{12}(u,v_{12})r^{13}(h,v_{13})+
[r^{12}(u,v_{12})-r^{12}(u-h,v_{12})]r^{13}(h,v_{13})=\\
&\frac{r^{23}(h,v_{23})-r^{23}(h,v_{13})}{v_{12}}cP^{12}+
[r^{23}(h,v_{23})\wt{r}^{12}(u,v_{12})-\wt{r}^{12}(u,v_{12})r^{13}(h,v_{13})]+\\
&[\wt{r}^{12}(u,v_{12})-\wt{r}^{12}(u-h,v_{12})]r^{13}(h,v_{13}).
\end{align*}
Passing to the limit $v_2\to v_1$ we derive
\begin{align*}
&r^{13}(u,v)r^{23}(-u+h,v)=-\frac{\pa r^{23}}{\pa v}(h,v)cP^{12}+[r^{23}(h,v)\wt{r}^{12}(u,0)-
\wt{r}^{12}(u,0)r^{13}(h,v)]+\\
&[\wt{r}^{12}(u,0)-\wt{r}^{12}(u-h,0)]r^{13}(h,v).
\end{align*}
Next, we are going to apply the operator $\mu\ot\id:A\ot A\ot A\to A\ot A$, where $\mu$ is
the product on $A$. We use the following easy observations:
$$(\mu\ot\id)(x^{13}y^{23})=xy, \ (\mu\ot\id)(x^{23}y^{12}-y^{12}x^{13})=0,\ 
(\mu\ot\id)(x^{23}P^{12})=1\ot \tr_1(x),$$
where $x,y\in A\ot A$, $\tr_1=\tr\ot\id:A\ot A\to A$ (the last property follows from the identity
$\sum_{ij} e_{ij}ae_{ji}=\tr(a)\cdot 1$ for $a\in A$).
Thus, applying $\mu\ot\id$ to the above equation we get
$$r(u,v)r(-u+h,v)=-c\cdot 1\ot\tr_1(\frac{\pa r}{\pa v}(h,v))+
(\mu\ot\id)\left([\wt{r}^{12}(u,0)-\wt{r}^{12}(u-h,0)]r^{13}(h,v)\right).$$
Finally, taking the limit $h\to 0$ we derive
\begin{equation}
s(u,v)=-c\cdot 1\ot\tr_1(\frac{dr_0(v)}{dv})+
\mu(\frac{\pa r(u,0)}{\pa u})\ot 1,
\end{equation}
where we used the equalities $\frac{\pa r(0,v)}{\pa v}=\frac{dr_0(v)}{dv}$ and
$\frac{\pa \wt{r}(u,v)}{\pa u}=\frac{\pa r(u,v)}{\pa u}$.
Hence, we can write $s(u,v)$ in the form
$$s(u,v)=a(u)\ot 1+1\ot b(v)+(f(u)+g(v))1\ot 1,$$
where $a(u)$ and $b(v)$ take values in $\splin_N$, and
$$b(v)=-c\pr\tr_1(\frac{dr_0(v)}{dv}).$$
The unitarity condition on $r(u,v)$ implies that $s^{21}(-u,-v)=s(u,v)$.
This immediately gives the required form of $s(u,v)$ with some $a\in\splin_N$, as well as the
formulas for $g(v)$, $f(u)$, $a$ and $\a(v)$.
The fact that $a$ is an infinitesimal symmetry of $r(u,v)$ follows from the second equality
in the identity of Lemma \ref{main-AYBE-lem1}.
\ed

\begin{lem}\label{last-AYBE-lem} 
Let $r_0(v)\in A\ot A$ be a unitary solution of the CYBE of the form
$$r_0(v)=\ov{r}_0(v)+\a(v)\ot 1-1\ot\a(-v)+h(v)\cdot 1\ot 1,$$
where $\ov{r}_0(v)\in\splin_N\ot\splin_N$ and $\a(v)\in\splin_N$.
Then 
$$[\ov{r}_0(v-v'),\a(v)\ot 1+1\ot \a(v')]=[\a(v),\a(v')]=0.$$ 
In particular, if $\a(v)$ depends
linearly on $v$, i.e., $\a(v)=b+v\cdot a$, then
$$[\ov{r}_0(v), a\ot 1+1\ot a]=[\ov{r}_0(v),b\ot 1+1\ot b+va\ot 1]=[b,a]=0.$$
\end{lem}

\Pf . Applying $\pr\ot\pr\ot\pr$ to both sides of the CYBE we see that $\ov{r}_0$ itself satisfies
the CYBE. Taking this into account the equation can be rewritten as
$$[\ov{r}_0^{12}(v_{12}),\a^1(v_{13})+\a^2(v_{23})]+[\a^1(v_{12}),\a^1(v_{13})]+
c.p.(1,2,3)=0,$$
where $v_{ij}=v_i-v_j$ (the omitted terms are obtained by cyclically permuting $1,2,3$). 
Applying the operator $\pr\ot\pr\ot\id$ gives 
$$[\ov{r}_0^{12}(v_{12}),\a^1(v_{13})+\a^2(v_{23})]=0.$$
Now returning to the above equality and applying $\pr\ot\id\ot\id$ we derive that $[\a(v),\a(v')]=0$.
\ed

\noindent
{\it Proof of Theorem \ref{QYBE-thm}.}
(i) In the case $N=1$ the statement is vacuous, so we can assume that $N>1$.
By Lemma \ref{pole-lem2}, $r_0(v)$ has a simple pole at $v=0$ with the residue $cP$, where
$c\in\C^*$. Projecting  to $\splin_N$ we deduce that $\ov{r}_0(v)$ is nondegenerate.

\noindent
(ii) By Lemma \ref{main-AYBE-lem1}, $r(u,v)$ satisfies the QYBE iff
$$s^{23}(u,v_{23})r^{13}(2u,v_{13})=r^{13}(2u,v_{13})s^{23}(-u,v_{23}).$$
Using the formula for $s(u,v)$ from Lemma \ref{main-AYBE-lem2} 
we see that this is equivalent to the equality
$$[r(u,v),1\ot a]=0$$
which is equivalent to $a=0$ by the nondegeneracy of $r(u,v)$. Note that by Lemma 
\ref{main-AYBE-lem2}, both conditions (b) and (c) are also equivalent to the equality $a=0$. 
It remains to show the equivalence of (d) with this equality. To this end we use the identity
\begin{equation}\label{r0-r1-eq}
r_0^{12}(v)r_0^{13}(v+v')-r_0^{23}(v')r_0^{12}(v)+r_0^{13}(v+v')r_0^{23}(v')=r_1^{12}(v)+r_1^{13}(v+v')+r_1^{23}(v')
\end{equation}
deduced by substituting the Laurent expansions in the first variable into \eqref{AYBE}.
Let us denote the expression in the left-hand-side of \eqref{r0-r1-eq} by $AYBE[r_0](v,v')$. Using
the relation between $r_0(v)$ and $\ov{r}_0(v)$ from Lemma \ref{main-AYBE-lem2} we obtain
\begin{align*}
&-cN\cdot(\pr\ot\pr\ot\pr)\left(AYBE[r_0](v,v')-AYBE[\ov{r}_0](v,v')\right)=\\
&v\ov{r}_0^{12}(v)a^3+v'\ov{r}_0^{23}(v')a^1+(v+v')\ov{r}_0^{13}(v+v')a^2,
\end{align*}
where $a^1=a\ot 1\ot 1$, etc. Note that \eqref{r0-r1-eq} implies that $(\pr\ot\pr\ot\pr)AYBE[r_0](v,v')=0$. 
Therefore, it suffices to prove that the equation
$$v\ov{r}_0^{12}(v)a^3+v'\ov{r}_0^{23}(v')a^1+(v+v')\ov{r}_0^{13}(v+v')a^2=0$$
on $a\in\splin_N$ implies that $a=0$.
Passing to the limit as $v\to 0$ and $v'\to 0$ we deduce from the above equality that
$$(\pr\ot\pr\ot\pr)[P^{12}a^3+P^{23}a^1+P^{13}a^2]=0.$$
Let $a=\sum a_{ij} e_{ij}$.
Looking at the coefficient with $e_{ij}\ot e_{ji}\ot e_{ij}$ we deduce that $a_{ij}=0$ for $i\neq j$.
Finally, looking at the projection to $e_{12}\ot e_{21}\ot\splin_N$ we deduce that
$a_{ii}$ does not depend on $i$, hence $a=0$.
\ed

\noindent
(iii) It suffices to prove that under our assumptions
the infinitesimal symmetry $a\in\splin_N$ appearing in Lemma \ref{main-AYBE-lem2} is equal
to zero. We only have to consider the case when $\ov{r}_0$ has a period, i.e.,
$\ov{r}_0(v+p)=\ov{r}_0(v)$ for some $p\in\C^*$. 
By Lemma \ref{last-AYBE-lem}, it remains to check that the equation
$$[\ov{r}_0(v),b\ot 1+1\ot b+va\ot 1]=0$$
on $a,b\in\splin_N$ implies that $a=0$. 
From the periodicity of $\ov{r}_0$ we derive that
$$[\ov{r}_0(v), a\ot 1]=0.$$
By the nondegeneracy of $\ov{r}_0$, it follows that $a=0$.
\ed

\section{Solutions of the AYBE associated with simple vector bundles on degenerations of elliptic
curves}
\label{bundle-sec}

Now let us review how solutions of the AYBE arise from geometric structures on elliptic curves and
their degenerations. Let $X$ be a nodal projective curve over $\C$ of
arithmetic genus $1$ such that the dualizing sheaf on $X$ is isomorphic to $\OO_X$. 
Let us fix such an isomorphism. Recall that a vector bundle $V$ on $X$ is called {\it simple} 
if $\End(V)=\C$.
The following result follows from Theorems 1 and 4 of \cite{P-AYBE}.

\begin{thm}\label{bundlethm}
Let $V_1$, $V_2$ be a pair of simple vector bundles on $X$ such that
$\Hom^0(V_1,V_2)=\Ext^1(V_1,V_2)=0$. Let $y_1,y_2$ be a pair of distinct
smooth points of $X$. Consider the tensor
$$r^{V_1,V_2}_{y_1,y_2}\in 
\Hom(V_{1,y_1}^*,V_{2,y_1}^*)\otimes \Hom(V_{2,y_2}^*,V_{1,y_2}^*)$$ 
corresponding to the following composition
$$\Hom(V_{1,y_1},V_{2,y_1})\lrar{\Res_{y_1}^{-1}}
\Hom(V_1,V_2(y_1))\lrar{\ev_{y_2}}\Hom(V_{1,y_2},V_{2,y_2}),$$
where $V_{i,y}$ denotes the fiber of $V_i$ at a point $y\in X$, 
the map
$$\Res_{y}:\Hom(V_1,V_2(y))\wt{\ra}\Hom(V_{1,y},V_{2,y})$$
is obtained by taking the residue at $y$, and
the map $\ev_y$ is the evaluation at $y$.
Then for a triple of simple bundles $(V_1,V_2,V_3)$ such that each pair satisfies the above assumptions and for a triple of distinct points $(y_1,y_2,y_3)$ one has 
\begin{equation}\label{Rid}
(r^{V_3V_2}_{y_1y_2})^{12}(r^{V_1V_3}_{y_1y_3})^{13}-
(r^{V_1V_3}_{y_2y_3})^{23}(r^{V_1V_2}_{y_1y_2})^{12}+
(r^{V_1V_2}_{y_1y_3})^{13}(r^{V_2V_3}_{y_2y_3})^{23}=0
\end{equation}
in
$\Hom(V_{1,y_1}^*,V_{2,y_1}^*)\ot\Hom(V_{2,y_2}^*,V_{3,y_2}^*)\ot\Hom(V_{3,y_3}^*,V_{1,y_3}^*)$.
In addition the following unitarity condition holds:
\begin{equation}
(r^{V_1V_2}_{y_1y_2})^{21}=-r^{V_2V_1}_{y_2y_1}.
\end{equation}
\end{thm}

\begin{rem} The tensor $r^{V_1,V_2}_{y_1,y_2}$ in the above theorem is a certain triple
Massey product in the derived category of $X$, and the equation \eqref{Rid} follows from the
appropriate $A_{\infty}$-axiom (see \cite{P-AYBE}).
\end{rem}

We are going to apply the above theorem for bundles $V_i$ of the form $V_i=V\otimes L_i$,
where $V$ is a fixed simple vector bundle of rank $N$ on $X$
and $L_i$ are line bundles in $\Pic^0(X)$,
the neutral component of $\Pic(X)$. Also, we let points $y_i$ vary
in a connected component $X_0$ of $X$.
Uniformizations of $X_0\cap X^{reg}$ and of $\Pic^0(X)$ allow to describe $V_i$'s and $y_i$'s
by complex parameters.
Thus, using trivializations of the bundles $V_{i,y_j}^*$ we can view the tensor
$r^{V_1,V_2}_{y_1,y_2}$ in the above theorem as a function of complex variables
$r(u_1,u_2;v_1,v_2)\in A\ot A$, where $A=\Mat(N,\C)$, $u_i$ describes $V_i$, $v_j$ describes $y_j$.
Note that equation \eqref{Rid}
reduces to the AYBE in the case when $r$ depends only on the differences
of variables, i.e., $r(u_1,u_2;v_1,v_2)=r(u_1-u_2,v_1-v_2)$.

A different choice of trivializations of $V_{i,y_i}^*$ would lead to the tensor $\wt{r}(u_1,u_2,v_1,v_2)$
given by
$$\wt{r}(u_1,u_2;v_1,v_2)=
(\varphi(u_2,v_1)\ot \varphi(u_1,v_2))r(u_1,u_2;v_1,v_2)
(\varphi(u_1,v_1)\ot\varphi(u_2,v_2))^{-1}$$
where $\varphi(u,v)$ is a function with values in $\GL_N(\C)$.
We say that tensor functions $\wt{r}$ and $r$ related in this way are {\it equivalent}.
Note that the condition for functions to depend only on the differences $u_1-u_2$
and $v_1-v_2$ is not preserved under these equivalences in general. However,
if $(a,b)$ is a pair of commuting infinitesimal symmetries of $r(u_1-u_2,v_1-v_2)$ then
taking $\varphi(u,v)=\exp(ua+vb)$ we do get a tensor function $\wt{r}$ that depends only on
the differences, namely,
$$\wt{r}(u,v)=\exp[u(1\ot a)+v(b\ot 1)]r(u,v)\exp[-u(a\ot 1)-v(b\ot 1)]$$
(this kind of equivalence shows up in Theorem \ref{main-thm}(ii)).

Since we are interested in trigonometric solutions, we will be using the multiplicative
variables $x_i=\exp(u_i)$, $y_i=\exp(v_i)$. The solutions of \eqref{Rid} that we are going to
construct in the next section
will be equivalent to those depending only on the differences $u_1-u_2$, $v_1-v_2$.
It will be convenient for us also to work with the intermediate form of the AYBE
\begin{equation}\label{AYBE2}
r^{12}((x')^{-1};y_1,y_2)r^{13}(xx';y_1,y_3)-
r^{23}(xx';y_2,y_3)r^{12}(x;y_1,y_2)+
r^{13}(x;y_1,y_3)r^{23}(x';y_2,y_3)=0
\end{equation}
for the tensor $r(x;y_1,y_2)\in A\ot A$, obtained from \eqref{Rid} in the case when
$r(x_1,x_2;y_1,y_2)=r(x_1/x_2;y_1,y_2)$.
The corresponding unitarity condition has form 
\begin{equation}\label{skew-eq}
r^{21}(x;y_1,y_2)=-r(x^{-1};y_2,y_1).
\end{equation}

\section{Simple vector bundles on cycles of projective lines}
\label{cycles-sec}

Let $X=X_0\cup X_1\cup\ldots\cup X_{n-1}$ be the union of $n$ copies of 
$\P^1$'s glued (transversally) in a configuration of type $\wt{A}_{n-1}$, 
so that the point $\infty$ on $X_j$ is identified
with the point $0$ on $X_{j+1}$ for $j=0,\ldots,n-1$ (where we identify indices with elements
of $\Z/n\Z$).
A vector bundle $V$ of rank $N$ on $X$ is given by a collection of vector bundles $V_j$ of rank $N$ on
$X_j$ along with isomorphisms $(V_j)_{\infty}\simeq (V_{j+1})_0$.
Since every vector bundle on $\P^1$ splits into a direct sum of line bundles, we can assume that
$$V_j=\OO_{\P^1}(m^j_1)\oplus\ldots\oplus\OO_{\P^1}(m^j_N)$$
for every $j=0,\ldots,n-1$. Thus, the splitting types are described by  
the $N\times n$-matrix of integers $(m^j_i)$. 

Let $(z_0:z_1)$ denote the homogeneous coordinates on $\P^1$. We will use the standard trivialization 
of the fiber of $\OO_{\P^1}(1)$ at $0=(1:0)\in\P^1$ (resp., at $\infty$) given by the generating section $z_0$ (resp., $z_1$). Note that a section $s\in\OO_{\P^1}(1)$ is uniquely determined by its
values $s(0)$ and $s(\infty)$ (namely, $s=s(0)z_0+s(\infty)z_1$).

Let us fix a splitting type matrix $\bm=(m^j_i)$.
For every $\la\in \C^*$ we define the rank-$N$ bundle $V^{\la}=V^{\la}(\bm)$ on $X$ by using standard trivializations
of $V_j=\oplus_{i=1}^n \OO(m^j_i)$ at $0$ and $\infty$ and setting the transition isomorphisms
$(V_j)_{\infty}\simeq (V_{j+1})_0$ to be identical for $j=0,\ldots, n-2$, and the last transition map to
be
$$\la C^{-1}:(V_0)_0\to (V_{n-1})_{\infty}$$
where $C$ is the cyclic permutation matrix: $Ce_i=e_{i-1}$, where we identify the set of indices with 
$\Z/NZ$. Note that in this definition only the cyclic order on the indices $\{1,\ldots,N\}$ is used.
In particular, if we cyclically permute the rows of the matrix $(m^j_i)$ (by replacing
$m^j_i$ with $m^j_{i+1}$) then we get the same vector bundle.

Lemma \ref{simple-crit-lem} below provides a criterion for simplicity of $V^{\la}(\bm)$.
This result is well-known (see \cite{BDG}, Theorem 5.3). For completeness we include the proof.
It is also known that every simple vector bundle on $X$ is isomorphic to some $V^{\la}(\bm)$ 
(see {\it loc. cit.}).
It will be convenient to extend 
the $N\times n$-matrix $(m^j_i)$ to the matrix with columns numbered by $j\in\Z$
using the rule $m^{j+n}_i=m^j_{i-1}$. 

\begin{lem}\label{simple-crit-lem}
The vector bundle $V^{\la}(\bm)$ is simple iff 
the following two conditions are satisfied:

\noindent
(a) {\it  the differences $m^j_i-m^j_{i'}$ for $i,i'\in\Z/N\Z$ take values only $\{-1, 0, 1\}$};

\noindent
(b) for every $i,i'\in\Z/N\Z$, $i\neq i'$, 
the $nr$-periodic infinite sequence 
$$(m^j_i-m^j_{i'}), \ \ j\in\Z$$
is not identically $0$, and
the occurrences of $1$ and $-1$ in it alternate.

Furthermore, if (a) and (b) hold then $V^{\la_1}(\bm)\simeq V^{\la_2}(\bm)$ iff $(\la_1/\la_2)^N=1$.
\end{lem} 

\Pf . First, we observe that 
if $m^j_i-m^j_{i'}=2$ then there exists 
a nonzero morphism $\OO_{\P^1}(m^j_{i'})\to\OO_{\P^1}(m^j_i)$ vanishing at $0$ and $\infty$.
Viewing it as an endomorphism of $V_j$ we obtain a non-scalar endomorphism of $V^{\la}$.
Hence, the condition (a) is necessary. From now on let us assume that (a) is satisfied.

A morphism $V^{\la_1}\to V^{\la_2}$ is given by a 
collection of morphisms $A_j:V_j\to V_j$, $j=0,\ldots,n-1$,
such that $A_j(\infty)=A_{j+1}(0)$ for $j=0,\ldots,n-2$ and
$$A_0(0)=\frac{\la_1}{\la_2} C A_{n-1}(\infty)C^{-1}.$$
We can write these maps as matrices
$A_j=(a_{ii'}^j)_{1\le i,i'\le N}$, where $a_{ii'}^j\in H^0(\P^1,\OO(m^j_i-m^j_{i'}))$.
Let us allow the index $j$ to take all integer values by using the rule
$a_{ii'}^{j+n}=a_{i-1,i'-1}^j$. Note that we still have $a_{ii'}^j\in H^0(\P^1,\OO(m^j_i-m^j_{i'}))$
because of our convention on $m^j_i$ for $j\in\Z$.
Then the equations on $(A_j)$ can be rewritten as
\begin{equation}\label{hom-system-eq}
a^j_{ii'}(0)=x^{\delta(j)}a^{j-1}_{ii'}(\infty)
\end{equation}
for all $i,i'\in\Z/N\Z$ and $j\in\Z$,
where $x=\la_1/\la_2$, and $\delta(j)=1$ for $j\equiv 0(n)$, $\delta(j)=0$ otherwise. 
Due to condition (a) we have
the following possibilities for each $a_{ii'}^j$:

\noindent
(i) if $m^j_i<m^j_{i'}$ then $a_{ii'}^j=0$;

\noindent
(ii) if $m^j_i=m^j_{i'}$ then $a_{ii'}^j$ is a constant, so $a_{ii'}^j(0)=a_{ii'}^j(\infty)$;

\noindent
(iii) if $m^j_i>m^j_{i'}$ then $a_{ii'}^j$ is a section of $\OO(1)$, so it is uniquely determined
by its values at $0$ and $\infty$, and these values can be arbitrary. 

From this we can immediately derive that (b) is necessary for $V^{\la}$ to be simple.
Indeed, if for some $i\neq i'$ we have $m^j_i=m^j_{i'}$ for all $j\in\Z$ then we can get a solution
of \eqref{hom-system-eq} with $x=1$ by setting $a^j_{i+k,i'+k}=1$ for all $j,k\in\Z$ and letting
the remaining entries to be zero. This would give a non-scalar endomorphism of $V^{\la}$. 
Similarly, if for some $i\neq i'$ and some segment $[j,k]\sub\Z$ we have
$$(m^j_i-m^j_{i'}, m^{j+1}_i-m^{j+1}_{i'},\ldots, m^k_i-m^k_{i'})=(1,0,\ldots,0,1)$$
then we get a solution of \eqref{hom-system-eq} with $x=1$ by setting 
$$a_{ii'}^j=z_1, a_{ii'}^{j+1}=1,\ldots, a_{ii'}^{k-1}=1, a_{ii'}^k=z_0$$
and letting the remaining entries to be zero.

Conversely, assume (a) and (b) hold. 
Then one can easily derive that $V^{\la}$ is simple by analyzing the system
\eqref{hom-system-eq} (with $x=1$). Indeed, 
let us show first that $a_{ii'}^j=0$ for $i\neq i'$. It follows from (b) that in the case $m^j_i=m^j_{i'}$
we can either find a segment $[j_1,j]\sub\Z$ such that $m^k_i=m^k_{i'}$ for $j_1<k<j$ and
$m^k_{j_1}<m^k_{j_1}$, or a segment $[j,j_2]\sub\Z$ such that $m^k_i=m^k_{i'}$ for $j<k<j_2$ and
$m^k_{j_2}<m^k_{j_2}$. In either case applying iteratively \eqref{hom-system-eq} we derive
that $a_{ii'}^j=0$ (recall that in this case $a_{ii'}^j$ is a constant). In the case $m^j_i>m^j_{i'}$ we
can find both segments $[j_1,j]$ and $[j,j_2]$ as above, so that \eqref{hom-system-eq} implies
that $a_{ii'}^j(0)=a_{ii'}^j(\infty)=0$. Hence, $a_{ii'}^j=0$. 
The remaining part of the system \eqref{hom-system-eq} shows that all $a_{ii}^j$ are equal
to the same constant, i.e., $V^{\la}$ has no non-scalar endomorphisms.

The above argument also shows that a morphism $(a_{ii'}^j):V^{\la_1}(\bm)\to V^{\la_2}(\bm)$
has $a_{ii'}^j=0$ (assuming conditions (a) and (b) hold), while the remaining components
$a_{ii}^j\in\C$ satisfy the equations
$$a_{ii}^j=x^{\de(j)}a_{ii}^{j-1}, \ \ a_{ii}^{j+n}=a_{i-1,i-1}^j,$$
where $x=\la_1/\la_2$. This system has a nonzero solution iff $x^N=1$, in which
case the solution gives an isomorphism $V^{\la_1}(\bm)\simeq V^{\la_2}(\bm)$.
\ed

\section{Computation of the associative $r$-matrix arising as a Massey product}
\label{comp-sec}

Henceforward, we always assume that the matrix $(m^j_i)$ satisfies the conditions
of Lemma \ref{simple-crit-lem}.
Given a pair of parameters $\la_1,\la_2\in \C^*$ and a pair of points $y,y'\in X_0\setminus\{0,\infty\}$
we want to describe explicitly the maps
$$\Res_y:\Hom(V^{\la_1},V^{\la_2}(y))\to \Hom(V^{\la_1}_y,V^{\la_2}_y),$$
$$\ev_{y'}:\Hom(V^{\la_1},V^{\la_2}(y))\to \Hom(V^{\la_1}_{y'},V^{\la_2}_{y'})$$
and especially the composition $\ev_{y'}\circ\Res_y^{-1}$ (for generic $\la_1,\la_2$).
We will identify the target spaces of both maps with $N\times N$-matrices using
trivializations of the relevant line bundles over $y$ induced by the appropriate power of 
$z_0\in H^0(\P^1,\OO(1))$. We also use the global $1$-form trivializing $\om_X$ that
restricts to $dz/z$ on each $\P^1\setminus\{0,\infty\}$ (where $z=z_1/z_0$).

A morphism $V^{\la_1}\to V^{\la_2}(y)$ is given by a 
collection of morphisms 
$$A_0: V_0\to V_0(y), A_1:V_1\to V_1,\ldots, A_{n-1}:V_{n-1}\to V_{n-1}$$
with same equations as before. 
Writing these maps as matrices we can view $\Hom(V^{\la_1},V^{\la_2}(y))$ as the
space of solutions of \eqref{hom-system-eq}, where  
$a_{ii'}^j\in H^0(\P^1,\OO(m^j_i-m^j_{i'}))$ for $j\not\equiv 0(n)$
and $a_{ii'}^j\in H^0(\P^1,\OO(m^j_i-m^j_{i'})(y))$ for $j\equiv 0(n)$. 

Since the component $X_0$ plays a special role, we will use a shorthand notation $m_i:=m^0_i$,
$a_{ii'}:=a_{ii'}^0$. Let us also set $b_{ii'}=\Res_y(a_{ii'})$. 
Recall that for every pair $i,i'\in\Z/N\Z$ we have the
following three possibilities.

\noindent
(i) If $m_i<m_{i'}$ then we have $a_{ii'}=\frac{yb_{ii'}}{z_1-yz_0}$, so that
\begin{equation}\label{hom-res-eq1}
a_{ii'}(0)=-b_{ii'}, a_{ii'}(\infty)=yb_{ii'}.
\end{equation}

\noindent
(ii) If $m_i=m_{i'}$ then $a_{ii'}=\frac{a_{ii'}(\infty)z-a_{ii'}(0)y}{z-y}$ (where $z=z_1/z_0$), so 
we get the relation
\begin{equation}\label{hom-res-eq2}
a_{ii'}(\infty)-a_{ii'}(0)=b_{ii'}.
\end{equation}

\noindent
(iii) If $m_i>m_{i'}$ then $a_{ii'}$ is uniquely determined by $a_{ii'}(0)$, $a_{ii'}(\infty)$ and $b_{ii'}$.
Namely, one can easily check that
$$a_{ii'}=\frac{z(b_{ii'}+a_{ii'}(0)-ya_{ii'}(\infty))-ya_{ii'}(0)}{z-y}\cdot z_0+
\frac{z a_{ii'}(\infty)}{z-y}\cdot z_1.
$$

Note that in the above three cases we also have the following expressions for
$a_{ii'}(y')$:
\begin{equation}\label{eval-eq}
a_{ii'}(y')=\cases \frac{yb_{ii'}}{y'-y}, & m_i<m_{i'},\\
\frac{y'b_{ii'}}{y'-y}+a_{ii'}(0)=\frac{yb_{ii'}}{y'-y}+a_{ii'}(\infty), & m_i=m_{i'},\\
\frac{y'b_{ii'}}{y'-y}+a_{ii'}(0)+y'a_{ii'}(\infty), & m_i>m_{i'}.\endcases
\end{equation}

To compute $\ev_{y'}\circ\Res_y^{-1}$ means to express all the entries $a_{ii'}(y')$ in terms of
$(b_{ii'})$. The above formula gives such an expression in the case $m_i<m_{i'}$; in the case 
$m_i=m_{i'}$ we need to know either $a_{ii'}(0)$ or
$a_{ii'}(\infty)$; and in the case $m_i>m_{i'}$ we need to know both. Of course, in the latter two
cases one has to use equations \eqref{hom-system-eq}. Then condition (b)
of Lemma \ref{simple-crit-lem} 
will guarantee that we get a closed formula for $a_{ii'}(y')$ in terms of all the
entries $b_{ii'}$. To organize the computation it is convenient to use the 
complete order on the set of indices $\{1,\ldots,N\}$ given by
\vspace{2mm}

\noindent $(\star)$
{\it $i\prec i'$ if either $m_i<m_{i'}$ or $m_i=m_{i'}$ and the first nonzero term in the sequence
$(m^j_i-m^j_{i'})$, $j=0,1,\ldots$, is negative}.
\vspace{2mm}

The fact that this is a complete order 
follows immediately from condition (b) of Lemma \ref{simple-crit-lem}. 
We will write $(ii')>0$ if $i\prec i'$ and $(ii')<0$ if $i\succ i'$. We will also use
the notation $-(i,i')=(i',i)$.

Let us define a partially defined operation on pairs of distinct indices in $\Z/N\Z$ 
by setting 
$$\tau(ii')=(i-1,i'-1)\text{ if }(i-1)\prec (i'-1)\text{ and }m^j_{i}=m^j_{i'}\text{ for }0<j<n.$$ 
Note that $\tau$ is one-to-one. We denote by $\tau^{-1}$ the
(partially defined) inverse and by $\tau^k$ the iterated maps.
Condition (b) of Lemma \ref{simple-crit-lem} implies that for every pair of distinct indices $(ii')$ there exists $k>0$ such that $\tau^k$ is not defined on $(ii')$.

\noindent
{\bf Case 1. Assume that $i\prec i'$, i.e., $(ii')>0$}. 
Then either $m_i<m_{i'}$, or there exists $j>0$ such that
$m^{j'}_i=m^{j'}_{i'}$ for $0\le j'<j$ and $m^j_i<m^j_{i'}$. In the first case we can use formula
\eqref{eval-eq}. In the second case we have
$a^{j'}_{ii'}=\const$ for $0<j'<j$, $j\not\equiv 0 (n)$, while $a^j_{ii'}=0$.
Therefore, using \eqref{hom-system-eq} and \eqref{hom-res-eq2} iteratively
we get the following expression for $a_{ii'}(\infty)$:
\begin{equation}\label{sum-eq1}
-a_{ii'}(\infty)=\sum_{k\ge 1}x^{-k}b_{\tau^k(ii')},
\end{equation}
where the summation is only over a finite number of $k$'s for which $\tau^k(ii')$ is defined.
This gives the following formula
\begin{equation}\label{answer-1-eq}
a_{ii'}(y')=\frac{yb_{ii'}}{y'-y}-\sum_{k\ge 1}x^{-k}b_{\tau^k(ii')}\ \ \text{  if }(ii')>0,
\end{equation}
that works also for the case $m_i<m_{i'}$ (since in this case $\tau$
is not defined on $(ii')$).

\noindent
{\bf Case 2. Assume that $i\succ i'$, i.e., $(ii')<0$}. 
Then either $m_i>m_{i'}$, or there exists $j<0$ such that
$m^{j'}_i=m^{j'}_{i'}$ for $j< j'\le 0$ and $m^j_i<m^j_{i'}$. 
Assume first that $m_i>m_{i'}$. Note that in this case there still exists $j<0$
with the above property, and in addition there is $k>0$ such that 
$m^{j'}_i=m^{j'}_{i'}$ for $0\le j'<k$ and $m^k_i<m^k_{i'}$ (by condition (b) of Lemma  \ref{simple-crit-lem}).
Using equations \eqref{hom-system-eq} and \eqref{hom-res-eq2} we derive
that \eqref{sum-eq1} still holds and also we have
\begin{equation}\label{sum-eq2}
a_{ii'}(0)=\sum_{k\ge 1}y^{\eps(\si\tau^{-k}\si(ii'))}x^kb_{\si\tau^{-k}\si(ii')},
\end{equation}
where $\si$ is the transposition: $\si(i,i')=(i',i)$, 
the summation is only over those $k$ for which $\tau^{-k}\si(ii')$ is defined, 
$\eps(ii')=1$ for $(ii')>0$ and $\eps(ii')=0$ otherwise.
This gives 
\begin{equation}\label{answer-2-eq}
a_{ii'}(y')=\frac{y'b_{ii'}}{y'-y}+\sum_{k\ge 1}y^{\eps(\si\tau^{-k}\si(ii'))}x^kb_{\si\tau^{-k}\si(ii')}-
y'\sum_{k\ge 1}x^{-k}b_{\tau^k(ii')}\ \ \text{ if }(ii')<0.
\end{equation}
We observe that this formula still works in the case $m_i=m_{i'}$ 
(the second summation becomes empty in this case).

\noindent
{\bf Case 3. Assume that $i=i'$}.
In this case we have relations
$$a_{ii}(0)=xa_{i+1,i+1}(0)+xb_{i+1,i+1}$$
for all $i\in\Z/N\Z$. Solving this linear system for $a_{ii}(0)$ we get
$$a_{ii}(0)=(1-x^N)^{-1}\sum_{k=1}^{N}x^kb_{i+k,i+k}.$$
Finally, we derive
\begin{equation}\label{answer-3-eq}
a_{ii}(y')=\frac{y}{y'-y}b_{ii}+(1-x^N)^{-1}\sum_{k=0}^{N-1}x^kb_{i+k,i+k}.
\end{equation}

Formulas \eqref{answer-1-eq}, \eqref{answer-2-eq} and \eqref{answer-3-eq}
completely determine the map $\ev_{y'}\circ\Res_y^{-1}$, so we can compute 
the associative $r$-matrix corresponding to the family of simple vector bundles $V^{\la}$ on $X$:
$$r(x;y,y')=r_{\const}(x,y/y')+\sum_{\a>0,k\ge 1} [-x^{-k} e_{-\tau^k(\a)}\otimes e_{\a}+
y^{\eps(-\tau^{-k}(\a))}x^k e_{\tau^{-k}(\a)}\otimes e_{-\a}-y'x^{-k}e_{-\tau^k(-\a)}\otimes e_{-\a}],
$$
where
\begin{equation}\label{r-const-eq}
\begin{array}{l}
r_{\const}(x,z)=\frac{z}{1-z} \sum_{\a>0} e_{-\a}\otimes e_{\a}+
\frac{1}{1-z} \sum_{\a>0} e_{\a}\otimes e_{-\a}+\\
\frac{z}{1-z}\sum_i e_{ii}\otimes e_{ii}+
(1-x^N)^{-1}\sum_i\sum_{k=0}^{N-1} x^k e_{i+k,i+k}\otimes e_{ii}.
\end{array}
\end{equation}
In these formulas $i$ is an element of $\Z/NZ$, and
$\a$ denotes a pair of distinct indices in $\Z/N\Z$. 
By a simple rearrangement of terms we can rewrite $r(x;y,y')$ in the following way:
\begin{equation}\label{r-gen-eq}
\begin{array}{l}
r(x;y,y')=r_{\const}(x,y/y')+\\
\sum_{\a>0,k\ge 1}[x^k e_{\a}\otimes e_{-\tau^k(\a)}
-x^{-k}e_{-\tau^k(\a)}\otimes e_{\a}+
yx^k e_{-\a}\otimes e_{-\tau^k(-\a)}-y'x^{-k} e_{-\tau^k(-\a)}\otimes e_{-\a}].
\end{array}
\end{equation}
Recall that this is a solution of \eqref{AYBE2} with the unitarity condition \eqref{skew-eq}.

\begin{ex} Assume that $n>N$ and the only nonzero entries of $(m^j_i)$ are
$m_1^N=m_2^{N-1}=\ldots=m_{N-1}^1=1$. Then the domain of definition of $\tau$ is empty,
so in this case
we have $r(x;y,y')=r_{\const}(x,y/y')$. Hence, $r_{\const}(\exp(u),\exp(v))$ is a solution of the AYBE.
\end{ex}

Later we will show that $r(\exp(u);\exp(v_1),\exp(v_2))$ 
is equivalent to an $r$-matrix depending only on the
difference $v_1-v_2$ (see Lemma \ref{equiv-lem}), so that it gives a solution of the AYBE.

\section{Associative Belavin-Drinfeld triples associated with simple vector bundles}
\label{comb-sec}

The right-hand side of \eqref{r-gen-eq} depends only on the parameters $x,y,y'$ and on a certain combinatorial structure on the set $S=\{1,\ldots,N\}$. We are going to show that this structure consists
of an {\it associative BD-structure} as defined in the introduction together with a compatible complete
order (see below). Later we will show that one can get rid of the dependence on a complete order
by passing to an equivalent $r$-matrix (see Lemma \ref{equiv-lem}). However, for purposes of studying splitting types of simple vector bundles on cycles of projective lines the full combinatorial structure
described below may be useful.

\begin{defi} We say that a {\it complete order on a set $S$ is compatible with the cyclic order} given
by a cyclic permutation $C_0$ (or simply {\it compatible with $C_0$}) if $C_0$ takes
every non-maximal element to the next element in this order. In other words, if we identify $S$ with the segment of integers $[1,N]$ preserving the complete order then $C_0(i)=i+1$ (where the indices are identified with $\Z/N\Z$). In this case we set $\a_0=(s_{\max},s_{\min})\in\Ga_{C_0}$, where
$s_{\min}$ (resp., $s_{\max}$) is the minimal (resp., maximal) element of $S$.
A choice of a complete order on $S$ compatible with $C_0$ is equivalent to a choice of an element
$\a_0\in\Ga_{C_0}$. By an {\it associative BD-structure on a completely ordered set} $S$ we mean
an associative BD-structure $(C_0,C,\Ga_1,\Ga_2)$ on $S$ such that the complete order is
compatible with $C_0$.
\end{defi}
 
Note that a choice of an associative BD-structure on the completely ordered set $[1,N]$
such that $\a_0\not\in\Ga_1$ and $\a_0\not\in\Ga_2$, is equivalent to a choice of a
Belavin-Drinfeld triple in $A_{N-1}$ equipped with an {\it associative structure} as defined in \cite{Sch}.

We will need the following characterization of associative BD-structures 
on completely ordered sets such that $\a_0\not\in\Ga_2$.

\begin{lem}\label{BD-simple-lem} 
Let $(S,<)$ be a completely ordered finite set 
equipped with a transitive cyclic permutation $C:S\to S$. Then to give an {\it associative BD-structure} on 
$S$ with $\a_0\not\in\Ga_2$ is equivalent to giving
a pair of subsets $P_1$ and $P_2$ in the set of pairs of distinct elements of $S$,
such that $(C\times C)(P_1)=P_2$ and the following properties are satisfied:

\noindent
(a) For every $(s,s')\in P_2$ one has $s<s'$.

\noindent
(b) Assume that $s<s'<s''$. If $(s,s'')\in P_1$ 
then $(s,s'), (s',s'')\in P_1$. The same property holds for $P_2$.
Also, if $(s',s)\in P_1$ then $(s',s''), (s'',s)\in P_1$ (resp., if $(s'',s')\in P_1$ then $(s'',s), (s,s')\in P_1$).
\end{lem}

The proof is left for the reader. Let us observe only
that property (b) assures that $P_{\iota}$ is determined by $\Ga_{\iota}=P_{\iota}\cap\Ga_{C_0}$,
where $\iota=1,2$. 

Now let us check that in the setting of section \ref{comp-sec}
we do get a completely ordered set with an associative BD-structure.

\begin{lem}\label{BD-lem} 
Let $(m^j_i)$ be a $N\times n$-matrix satisfying conditions of Lemma \ref{simple-crit-lem}.
Equip the set $S=\{1,\ldots,N\}$ with the complete order $\prec$ given by $(\star)$ and 
the cyclic permutation $C(i)=i-1$. Also, let 
$$P_1=\{(ii')\ |\ m^j_i=m^j_{i'}\text{ for }0<j<n\text{ and }C(i)\prec C(i')\}.$$
Then these data define an associative BD-structure with $\a_0\not\in\Ga_2$.
\end{lem}

\Pf . We use Lemma \ref{BD-simple-lem}.
The only question is why property (b) holds. Let $i\prec i'\prec i''$.

Assume first that $(i,i'')\in P_2$. Then $m^j_{i+1}=m^j_{i''+1}$ for $j\in [1,n-1]$.
Suppose there exists $j\in [1,n-1]$ such that $m^j_{i+1}\neq m^j_{i'+1}$. Consider
the maximal such $j$. We have either $m^j_{i+1}<m^j_{i'+1}$ or $m^j_{i'+1}<m^j_{i''+1}$.
By condition (b) of Lemma \ref{simple-crit-lem}, the former assumption contradicts to
$i\prec i'$, while the latter contradicts to $i'\prec i''$. Hence, $m^j_{i+1}=m^j_{i'+1}=m^j_{i''+1}$
for all $j\in [1,n-1]$, so that $(i,i'), (i',i'')\in P_2$.

Assume that $(i,i'')\in P_1$. Then $m^j_i=m^j_{i''}$ for $j\in [1,n-1]$.
Furthermore, since $i\prec i''$ and $i-1\prec i''-1$ we should have $m^0_i=m^0_{i''}$ (by condition
(b) of Lemma \ref{simple-crit-lem}).
Suppose there exists $j\in [0,n-1]$ such that $m^j_i\neq m^j_{i'}$. Consider the minimal such $j$.
We have either $m^j_i>m^j_{i'}$ or $m^j_{i'}>m^j_{i''}$. But the former contradicts to $i\prec i'$,
and the latter contradicts to $i'\prec i''$. Therefore, $m^j_i=m^j_{i'}=m^j_{i''}$ for all $j\in [0,n-1]$,
so that $(i,i'), (i',i'')\in P_1$.

Finally, 
assume that $(i',i)\in P_1$ (resp., $(i'',i')\in P_1$). Then  $m^j_i=m^j_{i'}$ (resp., $m^j_{i'}=m^j_{i''}$)
for $j\in [1,n-1]$.
Also, since $i'\succ i$ and $i'-1\prec i-1$ (resp., $i''\succ i'$ and $i''-1\prec i'-1$), we necessarily have 
$m^0_i<m^0_{i'}$ (resp., $m^0_{i'}<m^0_{i''}$). Hence, $m^0_{i'}=m^0_{i''}$ (resp., $m^0_i=m^0_{i'}$).
Suppose there exists $j\in [1,n-1]$ such that
$m^j_{i'}\neq m^j_{i''}$ (resp., $m^j_i\neq m^j_{i'}$). 
Consider the minimal such $j$. Since $i'\prec i''$ (resp., $i\prec i'$), we have
$m^j_i=m^j_{i'}<m^j_{i''}$ (resp., $m^j_i<m^j_{i'}=m^j_{i''}$). 
But this contradicts to condition (b) of Lemma \ref{simple-crit-lem}
(applied to $i$ and $i''$). Therefore, $m^j_i=m^j_{i'}=m^j_{i''}$ for all $j\in [1,n-1]$.
Since $m^0_{i'}=m^0_{i''}$ (resp., $m^0_{i}=m^0_{i'}$), we have $i'-1\prec i''-1$ (resp., $i-1\prec i'-1$), and hence $(i',i'')\in P_1$ (resp., $(i,i')\in P_1$).
Also, $i''-1\prec i-1$ (by condition (b) of Lemma \ref{simple-crit-lem}), so that $(i'',i)\in P_1$.
\ed



We will need below the following two operations on associative BD-structures.

\begin{defi} For an associative BD-structure $(C_0,C,\Ga_1,\Ga_2)$ on a finite set $S$
we define 

\noindent (i) the {\it opposite associative BD-structure} to be $(C_0^{-1},C,\si(\Ga_1),\si(\Ga_2))$,
where $\si$ is the permutation of factors in 
$S\times S$ (note that $\si(\Ga_{C_0})=\Ga_{C_0^{-1}}$);

\noindent (ii) the {\it inverse associative BD-structure} to be $(C_0,C^{-1},\Ga_2,\Ga_1)$.
\end{defi}

Note that under passing to the opposite associative BD-structure each set $P_{\iota}$, $\iota=1,2$, gets replaced with $\si(P_{\iota})$.

\begin{thm}\label{geom-BD-thm}
An associative BD-structure on a completely ordered finite set $S$ is obtained by the construction of Lemma \ref{BD-lem} from some matrix $(m^j_i)$ 
(satisfying conditions of Lemma \ref{simple-crit-lem})
iff $\a_0\not\in\Ga_2$ and
$C=C_0^k$ for some $k\in\Z$ (relatively prime to $N=|S|$).
\end{thm}

\Pf . {\bf ``Only if"}.
Let us denote by $t_i=\sum_{j=0}^{n-1} m^j_i$, $i=1,\ldots,N$, the sums of entries in the rows
of the matrix $(m^j_i)$. Then we claim that for $i\prec i'$ one has
\begin{equation}\label{diff-rows-eq}
t_i-t_{i'}=\cases -1, &\text{ if } i-1\succ i'-1,\\ 0, \text{ otherwise. }\endcases
\end{equation}
Indeed, assume first that $m^j_i=m^j_{i'}$ for all $j\in [0,n-1]$. Then
$i-1\prec i'-1$ and $t_i=t_{i'}$, so the above equation holds. Next,
assume that $m^j_i\neq m^j_{i'}$ for some $j\in [0,n-1]$.
Then the first nonzero term in the sequence $(m^j_i-m^j_{i'})_{j\in [0,n-1]}$ is $-1$.
Since $-1$'s and $1$'s in this sequence alternate, we have $t_i-t_{i'}=0$ (resp., $t_i-t_{i'}=-1$)
iff the last nonzero term in the sequence $(m^j_i-m^j_{i'})_{j\in [0,n-1]}$ is $1$ (resp., $-1$). But this happens precisely when the first nonzero term in $(m^j_{i-1}-m^j_{i'-1})_{j\ge 0}$ is $-1$ (resp., $1$),
so \eqref{diff-rows-eq} follows.

Now assume that $i\prec i'\prec i''$. Then it follows from \eqref{diff-rows-eq} that either
$C(i)\prec C(i')\prec C(i'')$ or $C(i')\prec C(i'')\prec C(i)$ or $C(i'')\prec C(i)\prec C(i')$.
Since this holds for every triple $(i,i',i'')$, it is easy to deduce that $C=C_0^k$ for some $k\in\Z$.

\noindent {\bf ``If"}. 
First, note that the construction of the associative BD-structure on a completely ordered set $S$ given in Lemma \ref{BD-lem} can be rewritten as follows. Assume we are given a transitive
cyclic permutation $C$ of $S$
and a matrix $(m^j_s)$, where $j\in [0,n]$, $s\in S$. Then we can extend the range of the index
$j$ to $\Z$ using the rule $m^{j+n}_s=m^j_{C(s)}$. Assuming that condition (b) of Lemma
\ref{simple-crit-lem} holds for this extended matrix we can proceed to define the complete order
by $(\star)$ and the set $P_1$ as in Lemma \ref{BD-lem}. Of course, we can always identify
$S$ with $\{1,\ldots,N\}$ in such a way that $C(i)=i-1$, so that we get to the setup of Lemma
\ref{BD-lem}. The advantage of the new point of view is that we can also consider
the set $S=\{1,\ldots,N\}$ with the cyclic permutation $C(i)=i-k$, where
$k\in\Z/N\Z$ is relatively prime to $N$. Then as was noted above we have to modify the definition
of the extended matrix by using the rule $m^{j+n}_i=m^j_{i-k}$.

Note that changing $(m^j_i)$ to $(-m^j_i)$ changes the associative BD-structure on $S$ 
to the opposite BD-structure, and the complete order on $S$ gets reversed.
Let us denote by $w_0:S\to S$ the permutation that reverses the order. Assume that we have
$C=C_0^k$. Then conjugating by $w_0$
the BD-structure associated with $(-m^j_i)$ we get a BD-structure that is obtained from the original
one by leaving the complete order the same, changing $C=C_0^k$ to $C_0^{-k}$, and replacing
$P_1$ with $(w_0\times w_0)\si(P_1)$. 
Therefore, it is enough to show that Lemma \ref{BD-lem} produces  
all associative BD-structures with $C=C_0^{-k}$, where $N/2\le k<N$.

Next, we describe a construction of a class of matrices $(m^j_i)$ 
satisfying conditions of Lemma \ref{simple-crit-lem}.
Fix $k$, relatively prime to $N$, such that $N/2\le k<N$.
Start with a sequence $(a_1,\ldots, a_N)$ such that $a_1=1$, $a_N=n-1$ (where $n>1$), 
and for every $i\in [1,N-1]$ one has either $a_{i+1}=a_i$ or $a_{i+1}=a_i+1$. 
Then set $m^0_i=1$ for $i\in [k+1,N]$, $m^{a_i}_{k+1-i}=1$ for $i=1,\ldots,N$,
and let the remaining entries to be zero. We are going to check that this matrix satisfies conditions of
Lemma \ref{simple-crit-lem} (with the modified definition of the extended matrix).

It is convenient to extend the range of the index $i$ to $\Z$ by the rule
$m^j_i=m^j_{i+N}$, so that we get a matrix $(m^j_i)$ with rows and columns numbered
by $\Z$. Let us consider the subset $\La\sub[k+1-N,N]\times [0,n-1]$ defined by
$$\La=([k+1,N]\times\{0\})\cup\{(k+1-i,a_i)\ |\ i=1,\ldots,N\}.$$
Then we have 
$$\{(i,j)\in\Z\times\Z\ |\ m^j_i\neq 0\}=\cup_{a\in\Z} \La_a,\text{ where}$$
$$\La_0=\cup_{b\in\Z} (\La+b(2N-k,n)),\ \ \La_a=\La_0+a(N,0).$$
Note that each $\La_a$ intersects each row once, and if we denote by $(i,j_a(i)))$ the 
intersection point of $\La_a$ with the $i$th row then either $j_a(i-1)=j_a(i)$ or $j_a(i-1)=j_a(i)+1$.
In other words, as we go down one row the point of intersection either stays in the same column or moves one step to the right.
It follows that the intersection of $\La_a$ with each column is a line segment. 
Moreover, it is easy to see that the number of elements in this intersection is at most $N$. Indeed, for columns corresponding to $j\equiv 0(n)$ the intersection segment has $N-k$ elements. 
On the other hand, for $j\not\equiv 0(n)$ this number 
is equal to the number of $i\in [1,N]$ such that $j\equiv a_i(n)$, so it is at most $N$.
This implies that $\La_a$ and $\La_{a'}$ are disjoint for $a\neq a'$.
Hence, $j_a(i)<j_{a+1}(i)$ for all $a\in\Z$ and $i\in\Z$.

Let us set $E_i=\{j_a(i)\ |\ a\in\Z\}$ for every $i\in\Z$. We have to check that for every pair of rows,
the $i$-th and the $i'$-th, where $i<i'<i+N$, one has $E_i\neq E_{i'}$, and the subsets
$E_i\setminus E_{i'}$ and $E_{i'}\setminus E_i$ in $\Z$ alternate. 

To prove that $E_i\neq E_{i'}$ we recall that by the construction, for every $b\in\Z$ 
the intersection of $\La_0$ with the $bn$-th column is the segment $[k+1+b(2N-k),N+b(2N-k)]$. 
The intersections of other sets $\La_a$ with the same column are obtained from
the above segment by shifts in $N\Z$. Since
$2N-k$ is relatively prime to $N$, it follows that for appropriate $b\in\Z$ the intersection of
$\cup_a \La_a$ with the $bn$-th column contains exactly one of the numbers $i$ and $i'$.
Hence, $bn$ belongs to exactly one of the sets $E_i$ and $E_{i'}$.

Finally, we have to prove that subsets $E_i\setminus E_{i'}$ and $E_{i'}\setminus E_i$
alternate. Note that for all $a$ we have $j_a(i')\le j_a(i)$.
Hence, our assertion would follow once we check that for every $a\in\Z$ one has $j_a(i)\le j_{a+1}(i')$. 
Suppose we have $j_{a+1}(i')<j_a(i)$. Then the intersection of $\La_{a+1}$ with the $j_a(i)$-th
column is a segment $[i_1,i_2]$, where $i<i_1\le i_2<i'$. Since $\La_{a+1}=\La_a+(N,0)$, the
intersection of $\La_a$ with the $j_a(i)$-th column is $[i_1-N,i_2-N]$. Hence, $i\le i_2-N<i'-N$, which
contradicts our assumptions on $i$ and $i'$.
 
Now given a BD-structure on a set $S=\{1,\ldots,N\}$ with the complete order $1<2<\ldots<N$
and the cyclic permutation $C=C_0^{-k}$ (where $N/2\le k<N$) we define the sequence
$(a_1,\ldots, a_N)$ as follows. Set $a_1=1$, and for $i=1,\ldots,N-1$ 
set 
$$a_{i+1}=\cases a_i & \text{ if } \a_{k-i}\in \Ga_1,\\ a_i+1 & \text{ otherwise},\endcases$$
where $\a_j=(j,j+1)$ (this uniquely defines $n$). It is easy to check that the corresponding matrix $(m^j_i)$
realizes our BD-structure. 
\ed

\section{Solutions of the AYBE and associative BD-structures}
\label{comb-sol-sec}

Let $(S,<,C,\Ga_1,\Ga_2)$ be a completely ordered finite set with an associative BD-structure
such that $\a_0\not\in\Ga_2$.
As in the introduction, for an element $\a=(i,j)\in S\times S$ we set $e_{\a}=e_{ij}\in A_S\simeq
\Mat_N(\C)$ (where $N=|S|$, the rows and columns are numbered by $S$).
We write $(i,j)>0$ (resp., $(i,j)<0$) if $i<j$ (resp., $i>j$).
Also, for $\a=(i,j)$ we set $-\a=(j,i)$.
Mimicking formulas \eqref{r-const-eq} and \eqref{r-gen-eq} we define
\begin{equation}\label{r-const-S-eq}
\begin{array}{l}
r_{\const}(x,z)=\frac{z}{1-z} \sum_{\a>0} e_{-\a}\otimes e_{\a}+
\frac{1}{1-z} \sum_{\a>0} e_{\a}\otimes e_{-\a}+\\
\frac{z}{1-z}\sum_{i\in S} e_{ii}\otimes e_{ii}+
(1-x^N)^{-1}\sum_{i\in S}\sum_{k=0}^{N-1} x^k e_{i,i}\otimes e_{C^k(i),C^k(i)},
\end{array}
\end{equation}
\begin{equation}\label{r-gen-S-eq}
\begin{array}{l}
r(x;y,y')=r_{\const}(x,y/y')+\\
\sum_{\a>0,k\ge 1}[x^k e_{\a}\otimes e_{-\tau^k(\a)}
-x^{-k}e_{-\tau^k(\a)}\otimes e_{\a}]+
\sum_{\a<0,k\ge 1}[yx^k e_{\a}\otimes e_{-\tau^k(\a)}-y'x^{-k} e_{-\tau^k(\a)}\otimes e_{\a}].
\end{array}
\end{equation}
In the last formula we use the operation $\tau$ defined on $P_1\sub S\times S$;
the summation is extended only over those $(k,\a)$ for which $\tau^k(\a)$ is defined. 
Below we will show that $r(x;y,y')$ is a solution of \eqref{AYBE2} (see Theorem \ref{comb-thm}). 
To deduce from this Theorem \ref{main-thm}(i) we will use the following simple observation.

\begin{lem}\label{equiv-lem} 
In the above situation the $A_S\ot A_S$-valued function
$$-r(\exp(\frac{u_1-u_2}{N});\exp(v_1),\exp(v_2))$$ 
is equivalent to the one given in Theorem \ref{main-thm}(i) for the inverse associative BD-structure
$(C_0,C^{-1},\Ga_2,\Ga_1)$,
where $u=u_1-u_2$ and $v=v_1-v_2$.
\end{lem}

\Pf . We can assume that $S=[1,N]$ (the segment of natural numbers) with the standard order. 
Let us set 
$$\varphi(v)e_j=\exp(-\frac{jv}{N})e_j.$$ 
Then the corresponding equivalent matrix $\wt{r}(u_1,u_2;v_1,v_2)$ is obtained from 
$r(\exp(\frac{u_1-u_2}{N});\exp(v_1),\exp(v_2))$ by multiplying
each term $e_{ij}\ot e_{j'i'}$ with $\exp(\frac{(j-i)v_1-(j'-i')v_2}{N})$.
Now we observe that $r_{const}(x,y/y')$ is a linear combination of 
$e_{ij}\ot e_{j'i'}$, where $j-i=j'-i'$. Such a term gets multiplied by $\exp(\frac{(j-i)(v_1-v_2)}{N})$.
The same is true about the terms in $r(x;y,y')$ not containing $y$ or $y'$.
Indeed, if $i<j$ and $\tau^k$ is defined on $(i,j)$ then $C^k(j)-C^k(i)=j-i$.
On the other hand, the terms involving $y=\exp(v_1)$ and $y'=\exp(v_2)$ are linear combinations of
$e_{i,j}\ot e_{j',i'}$, where $j'-i'=j-i+N$. Indeed, this follows from the fact
that if $i>j$ and $\tau^k$ is defined on $(i,j)$ then $C^k(j)-C^k(i)=j-i+N$ (the proof reduces to the
case $(i,j)=(N,1)$). The only other observation we use to rewrite $-\wt{r}$ in the form
given in Theorem \ref{main-thm}(i) (with $C$ replaced by $C^{-1}$ and $\Ga_1$ and $\Ga_2$ exchanged) is that for $0<m<N$ and for $i,j\in [1,N]$
we have $j-i\equiv m(N)$ iff either $i<j$ and $j=i+m$, or $i>j$ and $j=i+m-N$.
\ed

Since for every associative BD-structure on a finite set $S$ we can choose a compatible complete
order in such a way that $\a_0\not\in\Ga_2$, Theorem \ref{main-thm}(i) will follow easily from 
the above lemma and the next result.

\begin{thm}\label{comb-thm} 
Let $(S,<,C,\Ga_1,\Ga_2)$ be a completely ordered finite set with an associative BD-structure
such that $\a_0\not\in\Ga_2$. Then the function
$r(x;y,y')$ given by \eqref{r-gen-S-eq}
is a solution of \eqref{AYBE2} satisfying the unitarity condition 
\eqref{skew-eq}.
\end{thm}

\begin{rem}
By Theorem \ref{geom-BD-thm} we already know the statement to be true 
if $C=C_0^k$. Also, the work \cite{Sch} deals with the case when
in addition $\a_0\not\in\Ga_1$ (this fact will be used below).
\end{rem}

The rest of this section will be occupied with the proof of Theorem \ref{comb-thm} (in the
end we will also explain how to deduce Theorem \ref{main-thm}(i)).

Let us denote by $P=\sum_{i,j}e_{ij}\otimes e_{ji}$ the permutation tensor.
Then we can rewrite our $r$-matrix in the form
$$r(x;y,y')=a(x)+yb(x)-y'c(x)+\frac{y}{y'-y}P,$$
where
$$a(x)=(1-x^N)^{-1}\sum_{i\in S}\sum_{k=0}^{N-1} x^k e_{i,i}\otimes e_{C^k(i),C^k(i)}+
\sum_{\a>0} e_{\a}\otimes e_{-\a}+\sum_{\a>0,k\ge 1}[x^k e_{\a}\otimes e_{-\tau^k(\a)}
-x^{-k}e_{-\tau^k(\a)}\otimes e_{\a}],$$
$$b(x)=\sum_{\a<0,k\ge 1} x^k e_{\a}\otimes e_{-\tau^k(\a)},$$
$$c(x)=b^{21}(x^{-1})=\sum_{\a<0,k\ge 1} x^{-k}  e_{-\tau^k(\a)}\otimes e_{\a}.$$

\begin{lem}\label{Sch-lem} 
Assume that $\a_0\not\in\Ga_2$.
Let us set $\Ga'_1=\Ga_1\setminus\{\a_0\}$, $\Ga'_2=\tau(\Ga'_1)$. Then
$$a(x)+\frac{y}{y'-y}P$$
is exactly the $r$-matrix corresponding to the associative BD-structure
$(S,<,C,\Ga'_1)$. 
\end{lem}

\Pf . It is easy to see that $P'_1=\{\a\in P_1\ |\ \a>0\}$. Thus, the terms $b(x)$ and $c(x)$
in the $r$-matrix associated with the new associative BD-structure vanish.
We claim that the term $a(x)$ for the new associative BD-structure is the same as for
the old one. Indeed, it is enough to check that $\a\in P'_1$ is in the domain of definition of
$\tau^k$ iff it is in the domain of $(\tau')^k$, where $\tau':P'_1\to P'_2$ is the bijection
induced by $\tau$. But this follows immediately from the fact that $P_2$ consists
only of $\a>0$ (due to the assumption that $\a_0\not\in\Ga_2$).
\ed

Let us denote by $AYBE[r](x,x';y_1,y_2,y_3)$ the left-hand side of \eqref{AYBE2}.

\begin{lem}\label{AYBE-lem}
Consider the $r$-matrix of the form  
\begin{equation}\label{r-form-eq}
r(x;y_1,y_2)=a(x)+y_1b(x)-y_2c(x)+\frac{y_1}{y_2-y_1}P,
\end{equation}
where $a^{21}(x^{-1})+a(x)=P$ and $b^{21}(x^{-1})=c(x)$.
Then $r$ satisfies the unitarity condition \eqref{skew-eq}.
Also, $AYBE[r]=0$ iff the following equations are satisfied:

\noindent
(i) $AYBE[a]=0$;

\noindent
(ii) $b^{12}(x)b^{13}(x')=0$;

\noindent
(iii) $b^{13}(x)b^{23}(x')=b^{21}(x')b^{13}(xx')+b^{23}(xx')b^{12}(x)$;

\noindent
(iv) $c^{13}(x)a^{23}(x')+a^{12}((x')^{-1})c^{13}(xx')=c^{23}(xx')a^{12}(x)-a^{13}(x)c^{23}(x')$.
\end{lem}

\Pf . The unitarity 
condition follows immediately from our assumptions on $a(x)$, $b(x)$ and $c(x)$.
It is easy to check that
\begin{align*}
&AYBE[a(x)+y_1b(x)-y_2c(x)+\frac{y_1}{y_2-y_1}P](x,x';y_1,y_2,y_3)=\\
&AYBE[a(x)+y_1b(x)-y_2c(x)](x,x';y_1,y_2,y_3)
-y_1c^{21}(x')P^{13}-y_2c^{13}(x)P^{23}-y_1b^{23}(xx')P^{12}.
\end{align*}
Now the conditions (i)-(iv) are obtained by equating to zero coefficients
with various monomials in $y_1$, $y_2$ and $y_3$ (of degree $\le 2$).
Namely, (i) is obtained by looking at the constant term (i.e., by substituting $y_i=0$).
Conditions (ii), (iii) and (iv) are obtained by looking at the coefficients with $y_1^2$, $y_1y_2$ and
$y_3$, respectively. To see that these conditions imply $AYBE[r]=0$ we can use the identity
$$AYBE[r](x,x';y_2,y_3,y_1)^{231}=AYBE[r]((xx')^{-1},x;y_1,y_2,y_3)$$
that holds for any $r$ satisfying the unitarity condition \eqref{skew-eq}.
\ed

Let us introduce the following notation. For every $k\ge 1$
we denote by $P(k)\sub P_1$ the domain of definition of $\tau^k$ and by 
$P(k)^+\sub P(k)$ (resp., $P(k)^-\sub P(k)$)
the set of all $\a>0$ (resp., $\a<0$) contained in $P(k)$. Note that $P(1)=P_1$.
The assumption $\a_0\not\in\Ga_2$ implies that $\tau(P(k))\sub P(k-1)^+$.
Using this notation we can rewrite our formulas for $a(x)$, $b(x)$ and $c(x)$ as follows:
\begin{align*}
&a(x)=(1-x^r)^{-1}\sum_{0\le k<r,i} x^k e_{i,i}\otimes e_{C^k(i),C^k(i)}+
\sum_{i<j} e_{i,j}\otimes e_{j,i}+\\
&\sum_{(i,j)\in P(k)^+} [x^k e_{i,j}\otimes e_{C^k(j),C^k(i)}-x^{-k}e_{C^k(j),C^k(i)}\otimes e_{i,j}],
\end{align*}
$$b(x)=\sum_{k\ge 1, (i,j)\in P(k)^-}x^k e_{i,j}\otimes e_{C^k(j),C^k(i)}.$$
$$c(x)=\sum_{k\ge 1, (i,j)\in P(k)^-}x^{-k} e_{C^k(j),C^k(i)}\otimes e_{i,j}.$$

The following two combinatorial observations are also going to be useful in the proof.

\begin{lem}\label{comb-lem1} 
Let $(i_1,i_2,i_3)$ be a triple of elements of $S$ and let $k\ge 1$.
Then the following two conditions are equivalent:

\noindent
(a) $(i_1,i_3)\in P(k)^-$ and $i_1<i_2$ (resp., $i_2<i_3$);

\noindent
(b) $(i_1,i_2)\in P(k)^+$ and $(i_2,i_3)\in P(k)^-$ (resp., $(i_1,i_2)\in P(k)^-$ and $(i_2,i_3)\in P(k)^+$).
\end{lem}

The proof is straightforward and is left to the reader.

\begin{lem}\label{comb-lem2}
Let $k\ge 1$. Then for every $(i_1,i_2)\in P(k)^-$ one has a decomposition
$S=S_1\sqcup S_2$, where
$$S_1=\{i\ |\ i<i_1, C^k(i)>C^k(i_1)\},\ \ S_2=\{i\ |\ i>i_2, C^k(i)<C^k(i_2)\}.$$
\end{lem}

\Pf . We can assume that $S=[1,N]$ with the standard order.
Note that the map $C^k$ restricts to a bijection 
$$[i_1,N]\sqcup [1,i_2]\wt{\to}[C^k(i_1),C^k(i_2)].$$
Passing to the complements we derive that the open segment $(i_2,i_1)$ is the disjoint
union of its intersections with $S_1$ and $S_2$.
Next, if $i\le i_2$ then $(i_1,i)\in P(k)^-$ (by Lemma \ref{comb-lem1}), so that $C^k(i_1)<C^k(i)$. Hence, 
$[1,i_2]\sub S_1\setminus S_2$. Similarly, $[i_1,N]\sub S_2\setminus S_1$.
\ed

\noindent
{\it Proof of Theorem \ref{comb-thm}.}
Let us check that equations (i)-(iv) of Lemma \ref{AYBE-lem} hold in our case.
Equation (i) follows from Lemma \ref{Sch-lem} and Theorem 3.4 of \cite{Sch}. More precisely,
one can easily check that in the case when $\a_0\not\in\Ga_1$ and $\a_0\not\in\Ga_2$ our $r$-matrix coincides with the associative $r$-matrix constructed in \cite{Sch} for 
the opposite associative BD-structure on $S$.
Equation (ii) follows from the fact that for any $(i,j),(i',j')\in P(1)^-$ one has
$i'>j$ and $i>j'$ (otherwise we would have $\Ga_1=\Ga_S$).
To check equation (iii) we write
$$b^{13}(x)b^{23}(x')=
\sum_{k\ge 1, m\ge 1;  (i,j)\in P(k)^-, (i',j')\in P(m)^-; C^k(i)=C^m(j')}
x^k(x')^m e_{i,j}\otimes e_{i',j'}\otimes e_{C^k(j),C^m(i')}.$$ 
Note that we cannot have $k=m$ since this would imply that $i=j'$ contradicting the assumption that
$(i,j)\in P(k)^-\sub P(1)^-$ and $(i',j')\in P(m)^-\sub P(1)^-$.
Hence, we can split the summation into two parts: one with $k>m$ and one with $k<m$.
Denoting $k-m$ (resp., $m-k$) by $l$ in the first (resp., second) case, we can rewrite these sums as
$$\Si_1=\sum_{l\ge 1,m\ge 1; (i,j)\in P(m+l)^-, (i',C^l(i))\in P(m)^-}
x^l(xx')^m e_{i,j}\otimes e_{i',C^l(i)}\otimes e_{C^{m+l}(j),C^m(i')},$$
$$\Si_2=\sum_{l\ge 1,m\ge 1; (i',j')\in P(m+l)^-, (C^l(j'),j)\in P(m)^-}
(xx')^m(x')^l e_{C^l(j'),j}\otimes e_{i',j'}\otimes e_{C^m(j),C^{m+l}(i')}.$$
On the other hand, we have
$$b^{23}(xx')b^{12}(x)=
\sum_{l\ge 1, m\ge 1; (i,j)\in P(l)^-, (i',C^l(j))\in P(m)^-}
x^l(xx')^m e_{i,j}\otimes e_{i',C^l(i)}\otimes e_{C^{m+l}(j),C^m(i')}.$$
We claim that this is equal to $\Si_1$. Indeed, the condition $(i,j)\in P(m+l)^-$ is equivalent
to the conjuction of $(i,j)\in P(l)^-$ and $(C^l(i),C^l(j))\in P(m)^+$.
Now our claim follows from Lemma \ref{comb-lem1} applied to the triple $(i',C^l(i),C^l(j))$
(recall that $C^l(i)<C^l(j)$ since $(i,j)\in P(l)$).
Similarly, we check that $b^{21}(x')b^{13}(xx')=\Si_2$, which finishes the proof of equation (iii).

Finally, let us verify equation (iv). We can split both terms in the left-hand side 
of this equation into four sums
according to the four pieces comprising $a(x)$:
$$c^{13}(x)a^{23}(x')=L_1+L_2+L_3-L_4,\ \ a^{12}((x')^{-1})c^{13}(xx')=-L_5+L_6+L_7-L_8,$$
where
\begin{align*}
&L_1=(1-(x')^N)^{-1}\sum_{0\le m<N, k\ge 1;(i,j)\in P(k)^-}x^{-k}(x')^m 
e_{C^k(j),C^k(i)}\ot e_{C^{N-m}(j),C^{N-m}(j)}\ot e_{i,j},\\
&L_2=\sum_{m\ge 1;i<j,(i',j)\in P(m)^-}x^{-m}e_{C^m(j),C^m(i')}\ot e_{i,j}\ot e_{i',i},\\
&L_3=\sum_{k\ge 1,m\ge 1;(i,j)\in P(k)^+,(i',C^k(j))\in P(m)^-}x^{-m}(x')^k
e_{C^{k+m}(j),C^m(i')}\ot e_{i,j}\ot e_{i',C^k(i)},\\
&L_4=\sum_{k\ge 1,m\ge 1;(i,i')\in P(m)^-,(i',j)\in P(k)^+}x^{-m}(x')^{-k}
e_{C^m(i'),C^m(i)}\ot e_{C^k(j),C^k(i')}\ot e_{i,j},\\
&L_5=(1-(x')^N)^{-1}\sum_{0\le m<N,k\ge 1;(i,j)\in P(k)^-}(x')^{N-m}(xx')^{-k}
e_{C^k(j),C^k(i)}\ot e_{C^{k+m}(j),C^{k+m}(j)}\ot e_{i,j},\\
&L_6=\sum_{k\ge 1,(i,j)\in P(k)^-,i'<C^k(j)} (xx')^{-k}e_{i',C^k(i)}\ot e_{C^k(j),i'}\ot e_{i,j},\\
&L_7=\sum_{k\ge 1,m\ge 1;(i,j)\in P(k)^-,(i',C^k(j))\in P(m)^+}(x')^{-m}(xx')^{-k}
e_{i',C^k(i)}\ot e_{C^{k+m}(j),C^m(i')}\ot e_{i,j},\\
&L_8=\sum_{k\ge 1,m\ge 1;(i,j)\in P(k)^-,(i',j')\in P(m)^+,C^m(i')=C^k(j)}(x')^m(xx')^{-k}
e_{C^m(j'),C^k(i)}\ot e_{i',j'}\ot e_{i,j}.
\end{align*}
We split each of the sums $L_4$ and $L_8$ into $3$ parts according to the ranges of summation $k=m$, $k>m$, and $k<m$ 
(in the last two cases we make substitutions $k\mapsto k+m$ and $m\mapsto k+m$, respectively):
$$L_4=L_{4,1}+L_{4,2}+L_{4,3},\ \ L_8=L_{8,1}+L_{8,2}+L_{8,3},$$
where
\begin{align*}
&L_{4,1}=\sum_{k\ge 1,(i,j)\in P(k)^-,i'<j}(xx')^{-k}e_{C^k(i'),C^k(i)}\ot e_{C^k(j),C^k(i')}\ot e_{i,j},\\
&L_{4,2}=\sum_{k\ge 1,m\ge 1;(i,i')\in P(m)^-,(i',j)\in P(k+m)^+}
x^{-m}(x')^{-k-m}e_{C^m(i'),C^m(i)}\ot e_{C^{k+m}(j),C^{k+m}(i')}\ot e_{i,j},\\
&L_{4,3}=\sum_{k\ge 1,m\ge 1;(i,i')\in P(k+m)^-,(i',j)\in P(k)^+}
x^{-k-m}(x')^{-k}e_{C^{k+m}(i'),C^{k+m}(i)}\ot e_{C^k(j),C^k(i')}\ot e_{i,j},\\
&L_{8,1}=\sum_{k\ge 1,(i,j')\in P(k)^-,j<j'}x^{-k}e_{C^k(j'),C^k(i)}\ot e_{j,j'}\ot e_{i,j},\\
&L_{8,2}=\sum_{k\ge 1,m\ge 1;(i,j)\in P(k+m)^-,(C^k(j),j')\in P(m)^+}x^{-k-m}(x')^{-k}
e_{C^m(j'),C^{k+m}(i)}\ot e_{C^k(j),j'}\ot e_{i,j},\\
&L_{8,3}=\sum_{k\ge 1,m\ge 1;(i,C^m(i'))\in P(k)^-,(i',j')\in P(k+m)^+}x^{-k}(x')^m
e_{C^{k+m}(j'),C^k(i)}\ot e_{i',j'}\ot e_{i,C^m(i')}.
\end{align*}
Making appropriate substitutions of the summation variables and using Lemma \ref{comb-lem1} 
one can easily check that
$$L_2=L_{8,1}, \ L_3=L_{8,3}.$$
It follows that the left-hand side of (iv) is equal to
$$(L_1-L_5)+(L_6-L_{4,1})+(L_7-L_{4,2}) -(L_{4,3}+L_{8,2}).$$
Next, making the substitution $m\mapsto N-k-m$ in the sum $L_5$ we find 
$$L_1-L_5=-\sum_{0<m<k,(i,j)\in P(k)^-}x^{-k}(x')^{-m}
e_{C^k(j),C^k(i)}\ot e_{C^m(j),C^m(j)}\ot e_{i,j}.$$
Also, substituting $i'$ by $C^k(i')$ in $L_6$, switching $k$ and $m$ in $L_{4,2}$, and 
using Lemma \ref{comb-lem1} we find that
$$L_6-L_{4,1}=\sum_{k\ge 1,(i,j)\in P(k)^-,i'>j,C^k(i')<C^k(j)}(xx')^{-k}e_{C^k(i'),C^k(i)}\ot
e_{C^k(j),C^k(i')}\ot e_{i,j},$$
$$L_7-L_{4,2}=\sum_{k\ge 1,m\ge 1,(i,j)\in P(k)^-,(C^k(i'),C^k(j))\in P(m)^+,i'>j}
x^{-k}(x')^{-k-m}e_{C^k(i'),C^k(i)}\ot e_{C^{k+m}(j),C^{k+m}(i')}\ot e_{i,j}.$$
Finally, we can rewrite the sum of the other remaining terms as follows:
$$L_1-L_5-L_{4,3}-L_{8,2}=
-\sum_{k\ge 1,m\ge 1,(i,i',j)\in\Pi(k,m)}x^{-k-m}(x')^{-k}
e_{C^{k+m}(i'),C^{k+m}(i)}\ot e_{C^k(j),C^k(i')}\ot e_{i,j},$$
where $\Pi(k,m)$ is  the subset of $\{(i,i',j)\ |\ (i,j)\in P(k)^-, (C^k(i),C^k(i'))\in P(m)^+\}$
consisting of $(i,i',j)$ such that either $i'\le j$ or $C^k(j)<C^k(i')$.
It follows from Lemma \ref{comb-lem2} that 
$$\Pi(k,m)=\{(i,i',j)\ |\ (i,j)\in P(k)^-, (C^k(i),C^k(i'))\in P(m)^+,i'<i\}.$$

We deal similarly with the right-hand side of equation (iv). Namely, we write
$$c^{23}(xx')a^{12}(x)=R_1+R_2+R_3-R_4,\ \ a^{13}(x)c^{23}(x')=R_5+R_6+R_7-R_8,$$
where the parts correspond to the summands in $a(x)$. We also have a
decomposition $R_3=R_{3,1}+R_{3,2}+R_{3,3}$ (resp.,
$R_8=R_{8,1}+R_{8,2}+R_{8,3}$) obtained by collecting terms with 
$x^k(xx')^{-m}$ (resp., $x^{-k}(x')^{-m}$) with $k=m$, $k>m$ and $k<m$.
Now one can easily check that 
$$R_6=R_{3,1},\ R_7=R_{3,2}.$$
Also, we have
$$R_1-R_5=\sum_{m\ge 1,0<k\le m;(i,j)\in P(m)^-}x^{-k}(x')^{-m}e_{C^k(i),C^k(i)}\ot
e_{C^m(j),C^m(i)}\ot e_{i,j}.$$
We denote by $(R_1-R_5)_{k=m}$ and by $(R_1-R_5)_{k<m}$ parts of this sum corresponding
to the ranges $k=m$ and $k<m$.
Then we have
$$(R_1-R_5)_{k=m}+R_2+R_{8,1}=\sum_{k\ge 1, (i,j)\in P(k)^-, i\le i'\text{ or }C^k(i')<C^k(i)}
(xx')^{-k}e_{C^k(i'),C^k(i)}\ot e_{C^k(j),C^k(i')}\ot e_{i,j}.$$
Using Lemma \ref{comb-lem2} 
it is easy to see that the condition on $(i,j,i')$ in this summation can be replaced
by the conjuction of $(i,j)\in P(k)^-$, $j<i'$ and $C^k(i')<C^k(j)$ (same as in the formula for
$L_6-L_{4,1}$).
Finally, we have
\begin{align*}
&R_{8,2}-R_4=\\
&-\sum_{k\ge 1,m\ge 1,(i,j)\in P(m)^-,j'<i,(C^m(i),C^m(j'))\in P(k)^+}
x^{-k-m}(x')^{-m}e_{C^{k+m}(j'),C^{k+m}(i)}\ot e_{C^m(j),C^m(j')}\ot e_{i,j},
\end{align*}
$$(R_1-R_5)_{k<m}+R_{3,3}+R_{8,3}=\sum_{k\ge 1,m\ge 1,(i,j',j)\in\Pi'(k,m)}x^{-k}(x')^{-k-m}
e_{C^k(j'),C^k(i)}\ot e_{C^{k+m}(j),C^{k+m}(j')}\ot e_{i,j},$$
where $\Pi'(k,m)$ is the subset of $\{(i,j',j)\ |\ (i,j)\in P(k)^-,(C^k(j'),C^k(j))\in P(m)^+\}$ consisting
of $(i,j',j)$ such that either $i\le j'$ or $C^k(j')<C^k(i)$. By Lemma \ref{comb-lem2}, we get
$$\Pi'(k,m)=\{(i,j',j)\ |\ (i,j)\in P(k)^-,(C^k(j'),C^k(j))\in P(m)^+,j<j'\}.$$
Now it is easy to see that parts of the left-hand side and the right-hand side of equation (iv) match as follows:
\begin{align*}
L_6-L_{4,1} &=(R_1-R_5)_{k=m}+R_2+R_{8,1},\\
L_7-L_{4,2} &=(R_1-R_5)_{k<m}+R_{3,3}+R_{8,3},\\
L_1-L_5-L_{4,3}-L_{8,2} &=R_{8,2}-R_4.
\end{align*}
\ed

\noindent
{\it Proof of Theorem \ref{main-thm}(i).} As was already observed, the fact that $r(u,v)$ is a unitary
solution of the AYBE follows from Lemma \ref{equiv-lem} and Theorem \ref{comb-thm}.
It follows from Theorem \ref{QYBE-thm} that $r(u,v)$ also satisfies the QYBE for fixed $u$.
It remains to check the unitarity condition for the quantum $R$-matrix given by \eqref{Q-r-matrix}.
In view of the unitarity of $r(u,v)$ this boils down to proving the identity
\begin{equation}\label{s-eq}
s(u,v)=\left([\exp(\frac{v}{2})-\exp(-\frac{v}{2})]^{-2}-[\exp(\frac{u}{2})-\exp(-\frac{u}{2})]^{-2}\right)\cdot 
1\ot 1.
\end{equation}
To this end we observe that from Theorem \ref{QYBE-thm} and 
Lemma \ref{main-AYBE-lem2} we know that 
$$s(u,v)=(f(u)+g(v))\cdot 1\ot 1,$$
where $f(u)=\frac{1}{N}\tr\mu(\frac{\pa r(u,0)}{\pa u})$ and
$g(v)=-\frac{1}{N}(\tr\ot\tr)(\frac{\pa r(0,v)}{\pa v})$.
Now \eqref{s-eq} follows immediately from the equalities
$$f(u)=\frac{d}{du}\left(\frac{1}{\exp(u)-1}\right)=-[\exp(\frac{u}{2})-\exp(-\frac{u}{2})]^{-2},$$
$$g(v)=\frac{d}{dv}\left(\frac{1}{\exp(-v)-1}\right)=[\exp(\frac{v}{2})-\exp(-\frac{v}{2})]^{-2}.$$
\ed

\begin{rem} The following interesting observation is due to T.~Schedler. Assume that
$\Ga_1$ does not contain two consecutive elements of $\Ga_{C_0}$, say, $(C_0^{-1}(i_0),i_0)$
and $(i_0,C_0(i_0))$. Then the function $r(u,v)$ given by Theorem \ref{main-thm}(i) is equivalent
to the one of the form $\frac{1\ot 1}{\exp(u)-1}+{\bf r}(v)$. Indeed, let us denote by $O(i_0,i)$ the
minimal $k\ge 0$ such that $C^k(i_0)=i$. Then one can easily check that
$$a=\sum_i \frac{O(i_0,i)}{N} e_{ii}$$
is an infinitesimal symmetry of $r(u,v)$ and
$$\exp[u(1\ot a)]r(u,v)\exp[-u(a\ot 1)]=\frac{1\ot 1}{\exp(u)-1}+{\bf r}(v),$$
where ${\bf r}(v)$ depends only on $v$. Note that the fact that $r(u,v)$ is a unitary solution of the AYBE
is equivalent to the following equations on ${\bf r}(v)$:
$$AYBE[{\bf r}](v,v')={\bf r}^{13}(v+v'), \ \ {\bf r}^{21}(-v)+{\bf r}(v)=1\ot 1.$$
\end{rem}

\section{Meromorphic continuation}
\label{pole-sec}

As was shown in the proof of Theorem 6 of \cite{P-AYBE} (see also Lemma 4.14 of \cite{Sch}), 
a unitary solution of the
AYBE with the Laurent expansion \eqref{Laur-exp} at $u=0$ is uniquely determined by $r_0(v)$.
Therefore, it is not surprising that some of the results from \cite{BD} about solutions of the CYBE 
(such as meromorphic continuation) can be extended to solutions of the AYBE. 

First, we apply the above uniqueness principle to infinitesimal symmetries.

\begin{lem}\label{sym-lem} Let $r(u,v)$ be a nondegenerate unitary solution of the AYBE with
the Laurent expansion \eqref{Laur-exp} at $u=0$.
Then the algebras of infinitesimal symmetries of $r(u,v)$ and of $r_0(v)$ are the same (and
are contained in the algebra of infinitesimal symmetries of $\ov{r}_0$).
If in addition $\ov{r}_0$ has a period then these coincide with the commutative algebra of
infinitesimal symmetries of $\ov{r}_0$.
\end{lem}

\Pf . Let $a\in A$ be an infinitesimal symmetry of $r_0(v)$.
Then for any $t\in\C$ the function
$$\exp[t(a\ot 1+1\ot a)]r(u,v)\exp[-t(a\ot 1+1\ot a)]$$ 
is a solution of the AYBE with the same $r_0$-term in the Laurent expansion
at $u=0$. By the uniqueness mentioned above this implies that $\exp[t(a\ot 1+1\ot a)]$ 
commutes with $r(u,v)$, so $a$ is an infinitesimal symmetry of $r(u,v)$. 
Recall that by Theorem \ref{QYBE-thm}(i), $\ov{r}_0(v)$ is nondegenerate.
It is easy to see that if $\ov{r}_0$ is either elliptic or
trigonometric then the algebra of infinitesimal symmetries of $\ov{r}_0$ is commutative.
Indeed, in the elliptic case this algebra is trivial (see Lemma 5.1 of \cite{P-AYBE}). In the trigonometric case this follows from the fact proven in \cite{BD}
that there exists a pole $\ga$ of $\ov{r}_0$ such that
$$\ov{r}_0(v+\ga)=(\phi\ot\id)(\ov{r}_0(v)),$$
where $\phi$ is a Coxeter automorphism of $\splin_N$.
Thus, any infinitesimal symmetry is contained in the commutative algebra of $\phi$-invariant elements.
\ed

\begin{prop}\label{mer-cont-prop}
Assume $N>1$.
Let $r(u,v)$ be a nondegenerate unitary solution of the AYBE with the Laurent expansion 
\eqref{Laur-exp} at $u=0$, such that the equivalent conditions of Theorem \ref{QYBE-thm}(ii) hold. Then
$r(u,v)$ admits a meromorphic continuation to $D\times\C$, where $D$ is a neighborhood of
$0$ in $\C$. If $r(u,v)$ has a pole at $v=\ga$ then this pole is simple and
$\ov{r}_0(v)$ also has a pole at $v=\ga$.
\end{prop}

\Pf . Note that $\ov{r}_0(v)$ has a meromorphic continuation to $\C$ with at most simple poles
by Theorem 1.1 of \cite{BD}.
First, we want to deduce a meromorphic continuation for $r_0(v)$.
From the condition (c) in Theorem \ref{QYBE-thm}(ii) we know 
that 
\begin{equation}\label{r0-for}
r_0(v)=\ov{r}_0(v)+b\ot 1-1\ot b+h(v)\cdot 1\ot 1,
\end{equation} 
where $b$ is an infinitesimal
symmetry of $\ov{r}_0(v)$ (by Lemma \ref{last-AYBE-lem}). 
Note that $b$ is also an infinitesimal symmetry of $r_0(v)$. Hence, by Lemma \ref{sym-lem}, 
$b$ is an infinitesimal symmetry of $r(u,v)$. Applying the equivalence transformation
$$r(u,v)\mapsto  \exp[u(1\ot b)]r(u,v)\exp[-u(b\ot 1)]$$
we can assume that $b=0$. In this case we have 
$$AYBE[r_0](v_{12},v_{23})-AYBE[\ov{r}_0](v_{12},v_{23})
\equiv [h(v_{13})-h(v_{23})]\ov{r}_0^{12}(v_{12})+c.p.(1,2,3)\mod(\C\cdot 1\ot 1\ot 1)
$$
where we use the notation from the proof of Theorem \ref{QYBE-thm}
(the omitted terms are obtained by cyclically permuting $(1,2,3)$; we denote $v_{ij}=v_i-v_j$).
Applying $\pr\ot\pr\ot\id$ and using \eqref{r0-r1-eq} we obtain
\begin{equation}\label{mer-cont-eq}
[h(v+v')-h(v')]\ov{r}_0^{12}(v)=[(\pr\ot\pr)r_1(v)]^{12}-
(\pr\ot\pr\ot\id)AYBE[\ov{r}_0](v,v').
\end{equation}
Note that $AYBE[\ov{r}_0](v,v')$ is meromorphic on the entire $\C\times\C$ and has at most
simple poles at $v=\ga$, $v'=\ga$ and $v+v'=\ga$, where $\ga$ is a pole of $\ov{r}_0(v)$.
Also, by Lemma \ref{pole-lem2}, $r_1(v)$ is holomorphic near $v=0$.
Choose a small disk $D$ around zero such that
$r_1(v)$ is holomorphic in $D$ and $\ov{r}_0(v)$ has no poles or zeros in $D\setminus\{0\}$.
Assume that we already have a meromorphic continuation of $h(z)$ to some open subset $U\sub\C$
containing zero.
Then the above formula gives a meromorphic continuation of $h(z)$ to $U+D$. 
Iterating this process we continue $h(z)$ meromorphically to the entire complex plane.
Furthermore, it is clear from \eqref{mer-cont-eq} that $h(v)$ 
has only simple poles and is holomorphic outside the set of poles of $\ov{r}_0(v)$. 
Therefore, the same is true for $r_0(v)$.

Next, considering the constant terms of the Laurent expansions of the AYBE in $u'$ (keeping $u$ fixed)
we get
\begin{equation}\label{AYBE-der-eq}
r_0^{12}(v_{12})r^{13}(u,v_{13})+r^{13}(u,v_{13})r_0^{23}(v_{23})-r^{23}(u,v_{23})r^{12}(u,v_{12})=
\frac{\pa r^{13}}{\pa u}(u,v_{13}).
\end{equation}
Since we already know that $r_0(v)$ is meromorphic on the entire $\C$, we can use this equation
to get a meromorphic continuation of $r(u,v)$. Indeed, assume that $r(u,v)$ is meromorphic
in $D\times D$ for some open disk around zero $D\sub\C$. 
For fixed $v_{13}\in D$ 
the above equation gives a meromorphic extension of 
$$r^{23}(u,v_{21}+v_{13})r^{21}(-u,v_{21})=-r^{23}(u,v_{23})r^{12}(u,v_{12})$$
to $D\times\C$. By the nondegeneracy of $r(u,v)$ this allows to extend meromorphically
$r(u,v)$ from $D\times U$ to $D\times (U+D)$. Iterating this process we get the required
meromorphic extension. The assertion about poles follows easily from \eqref{AYBE-der-eq}
by fixing $v_{13}$ such that $r(u,v)$ has no pole at $v=v_{13}$ and 
$r(u,v_{13}-\ga)$ is nondegenerate, and considering the polar parts at $v_{12}=\ga$.
\ed

The argument in the following Lemma is parallel to that in Proposition 4.3 of \cite{BD}.

\begin{lem}\label{pole-aut-lem}
With the same assumptions as in Proposition \ref{mer-cont-prop}
for every pole $\ga$ of $\ov{r}_0(v)$ there exists an algebra
automorphism $\phi_{\ga}$ of $A$ and a constant $\la\in\C$ such that
$$r(u,v+\ga)=\exp(\la u)(\phi_{\ga}\ot\id)(r(u,v)).$$
\end{lem}

\Pf . From Proposition \ref{mer-cont-prop}
we know that the pole of $r(u,v)$ at $v=\ga$ is simple.
Set $\tau(u)=\lim_{v\to\ga}(v-\ga)r(u,v)$. Recall that $\lim_{v\to 0}vr(u,v)=cP$ for
$c\in\C^*$.
Let us define an operator $\phi(u)\in\End(A)$ by the equality
$$\tau(u)=(\phi(u)\ot\id)(cP).$$
Considering polar parts near $v=\ga$ in \eqref{AYBE} we get
$$\tau^{12}(-u')r^{13}(u+u',v'+\ga)=r^{23}(u+u',v')\tau^{12}(u).$$
The right-hand side can be rewritten as follows:
$$r^{23}(u+u',v')\tau^{12}(u)=c(\phi(u)\ot\id\ot\id)(r^{23}(u+u',v')P^{12})=
c(\phi(u)\ot\id\ot\id)(P^{12}r^{13}(u+u',v')),$$
Hence, we have
\begin{equation}\label{pole-comp-eq}
\tau^{12}(-u')r^{13}(u+u',v'+\ga)=c(\phi(u)\ot\id\ot\id)(P^{12}r^{13}(u+u',v')).
\end{equation}
Taking the residues at $v'=0$ we find
$$\tau^{12}(-u')\tau^{13}(u+u')=c^2(\phi(u)\ot\id\ot\id)(P^{12}P^{13}).$$
This means that $\phi(u)$ satisfies the identity 
$$\phi(u_1+u_2)(XY)=\phi(u_1)(X)\phi(u_2)(Y),$$
where $X,Y\in A$. Let $D$ be a small disk around zero in $\C$ such that
$\phi(u)$ is holomorphic on $D\setminus\{0\}$.
For every $u\in D\setminus\{0\}$ we denote by $I(u)\sub A$ the kernel of $\phi(u)$.
Then from the above identity we derive that $I(u)A\sub I(u+u')$ and $AI(u)\sub I(u+u')$
whenever $u,u',u+u'\in D\setminus\{0\}$. In particular, we deduce that $I(u)\sub I(u+u')$,
so $I(u)=I\sub A$ does not depend on $u\in D\setminus\{0\}$. It follows that $I$ is a two-sided ideal in 
$A$. 
Since $\phi(u)$ is not identically zero, we derive that $I=0$. Therefore, $\phi(u)$ is invertible for every
$u\in D\setminus\{0\}$. Now as in the proof of Lemma \ref{nondeg-no-pole-lem} we derive that
$$\phi(u)=\exp(\la u)\phi_{\ga}$$
for some $\la\in\C$, where $\phi_{\ga}$ is an algebra automorphism of $A$.
Applying $\phi^{-1}_{\ga}\ot\id\ot\id$ to \eqref{pole-comp-eq} we derive
$$\exp(-\la u')P^{12}(\phi^{-1}_{\ga}\ot\id\ot\id)r^{13}(u+u',v'+\ga)=\exp(\la u)P^{12}r^{13}(u+u',v').$$
This implies the required identity.
\ed

\begin{lem}\label{period-lem}
Keep the same assumptions as in Proposition \ref{mer-cont-prop}.
Assume that $\ov{r}_0(v+p)=\ov{r}_0(v)$ for some $p\in\C^*$. Then 
$r(u,v+p)=\exp(\la u)r(u,v)$ for some constant $\la\in\C$.
\end{lem}

\Pf . Consider the decomposition \eqref{r0-for} again.
The identity \eqref{mer-cont-eq} implies that $h(v+v')-h(v')$ is periodic in $v'$ with the period $p$.
Hence, $h(v+p)=h(v)+\la$ for some $\la\in\C$. It follows that 
$r_0(v+p)=r_0(v)+\la\cdot 1\ot 1$. Applying the rescaling $r(u,v)\mapsto\exp(-\la uv)r(u,v)$
we can assume that $r_0(v+p)=r_0(v)$. Now Lemma \ref{pole-aut-lem} implies that
$r(u,v+p)=(\phi_p\ot\id)r(u,v)$, where $\phi_p$ is an automorphism of $A$. Since
$r_0(v)$ is nondegenerate (as follows from Lemma \ref{pole-lem2}), we derive that $\phi_p=\id$.
\ed

We will use the following result in the proof of Theorem \ref{AYBE-v}. 

\begin{prop}\label{r0-prop} 
Assume $N>1$.
Let $r(u,v)$ be a nondegenerate unitary solution of the AYBE
with the Laurent expansion at $u=0$ of the form \eqref{Laur-exp}
such that the equivalent conditions of Theorem \ref{QYBE-thm}(ii) hold. 
Then one has
$$r_0(v)=\ov{r}_0(v)+b\ot 1+1\ot b+h(v)\cdot 1\ot 1,$$ 
$$h(v)=\la v+c h_0(c'v),$$ 
where $b\in\splin_N$ is an infinitesimal symmetry of $\ov{r}_0(v)$, $\la\in\C$, $c,c'\in\C^*$,
and $h_0(v)$ is one of the following three functions: Weierstrass zeta function $\zeta(v)$ associated with a lattice in $\C$; $\frac{1}{2}\coth(\frac{v}{2})$; or $\frac{1}{v}$. 
Furthermore, if $\ov{r}_0(v)$ is equivalent to a rational solution of the CYBE then $h_0(v)=\frac{1}{v}$.
\end{prop}

\Pf . The equation \eqref{r0-r1-eq} implies that
\begin{equation}\label{square-eq}
[r_0^{12}(v_{12})+r_0^{23}(v_{23})+r_0^{31}(v_{31})]^2=x^{12}(v_{12})+x^{23}(v_{23})+
x^{31}(v_{31}),
\end{equation}
where $x(v)=r_0(v)^2-2r_1(v)$ (and $v_{ij}=v_i-v_j$).
On the other hand, it is easy to see that $x(v)$ is the constant term of the Laurent expansion of 
$s(u,v)=r(u,v)r(-u,v)$ at $u=0$. 
Rescaling $r(u,v)$ we can assume that its residue at $v=0$ is equal to $P$ 
(see Lemma \ref{pole-lem2}).
Then we have
$$s(u,v)=[f(u)+g(v)]\cdot 1\ot 1,$$
where 
$$f(u)=\frac{1}{N}\tr\mu(\frac{\pa r(u,0)}{\pa u}),\ \ g(v)=-\frac{1}{N}(\tr\ot\tr)(\frac{dr_0(v)}{dv})$$
(see Lemma \ref{main-AYBE-lem2}).
If we change $r(u,v)$ to $\exp(\la uv)r(u,v)$ for some $\la\in\C$ then $f(u)$ changes to $f(u)+N\la$
(this operation also changes $r_0(v)$ to $r_0(v)+\la v\cdot 1\ot 1$).
Therefore, we can assume that $f(u)$ has no constant term in the Laurent expansion at $u=0$.
In this case we obtain $x(v)=g(v)\cdot 1\ot 1$.
Hence, denoting 
$$T(v_1,v_2,v_3)=r_0^{12}(v_{12})+r_0^{23}(v_{23})+r_0^{31}(v_{31})$$
we can rewrite \eqref{square-eq} as follows:
\begin{equation}\label{square-eq2}
T(v_1,v_2,v_3)^2=[g(v_{12})+g(v_{23})+g(v_{31})]\cdot 1\ot 1.
\end{equation}
Viewing $T(v_1,v_2,v_3)\in A\ot A\ot A$ as an endomorphism of $V\ot V\ot V$,
where $A=\End(V)$, we obtain
\begin{equation}\label{T-eq}
T(v_1,v_2,v_3)=T_0(v_1,v_2,v_3)+[h(v_{12})+h(v_{23})+h(v_{31})]\cdot\id_{V\ot V\ot V},
\end{equation}
where $T_0$ is a traceless endomorphism and $h(v)$ is defined from the decomposition
\eqref{r0-for}. Note also that for fixed (generic) $v_2$ and $v_3$ we have
$$\lim_{v_1\to v_2}(v_1-v_2)T(v_1,v_2,v_3)=P^{12}.$$
The latter operator has $S^2V\ot V$ and $\We^2 V\ot V$ as eigenspaces.
Therefore, for $v_1$ close to $v_2$ we have a decomposition
$$V\ot V\ot V=W_1\oplus W_2,$$
where $\dim W_1=N^2(N+1)/2$, $\dim W_2=N^2(N-1)/2$, and
$$(T(v_1,v_2,v_3)-\la\id)(W_1)=0,\ \ (T(v_1,v_2,v_3)+\la\id)(W_2)=0\text{, where }
\la^2=g(v_{12})+g(v_{23})+g(v_{31}).$$
Comparing the traces of both sides of \eqref{T-eq} we derive
$$\la=N[h(v_{12})+h(v_{23})+h(v_{31})].$$
Since $g(v)=-Nh'(v)$, we obtain
$$N[h(v_{12})+h(v_{23})+h(v_{31})]^2+h'(v_{12})+h'(v_{23})+h'(v_{31})=0.$$
Replacing $h(v)$ by $h(Nv)$ we get
\begin{equation}\label{h-eq}
[h(v_{12})+h(v_{23})+h(v_{31})]^2+h'(v_{12})+h'(v_{23})+h'(v_{31})=0.
\end{equation}
We are interested in solutions of this equation for an odd meromorphic function $h(v)$ in the neighborhood of zero having a simple pole at $v=0$. It is easy to see that the Laurent expansion of 
$h(v)$ at $v=0$ should have form $h(v)=1/v+c_3v^3+\ldots$. As shown in the proof of Theorem 5
of \cite{P-AYBE}, all such solutions of \eqref{h-eq} have form $c\cdot h_0(cv)$, where $h_0$ is one of the three functions described in the formulation.\footnote{Solutions of \eqref{h-eq} were
first described by L.~Carlitz in \cite{Car}.}

Finally, if $\ov{r}_0(v)$ is rational then its only pole is $v=0$ (see \cite{BD}). Therefore, by
Proposition \ref{mer-cont-prop}, $r_0(v)$ also cannot have poles outside zero, which implies
that $h_0(v)=\frac{1}{v}$.
\ed

\begin{rem} In the case when $\ov{r}_0(v)$ is either elliptic or trigonometric the assertion
of the above proposition can also be deduced from the explicit formulas for
$r(u,v)$ (the elliptic case is discussed in \cite{P-AYBE}, sec.2; the trigonometric case is considered
in Theorem \ref{main-thm}).
\end{rem}

\section{Classification of trigonometric solutions of the AYBE}
\label{class-sec}

Recall (see \cite{BD}) that to every nondegenerate trigonometric solution $\ov{r}_0(v)$
of the CYBE for $\splin_N$ with poles exactly at $2\pi i\Z$
Belavin and Drinfeld associate an automorphism of the Dynkin diagram $A_{N-1}$ by
considering the class of the automorphism $\phi$ of $\splin_N$ defined by
\begin{equation}\label{trig-aut-pole-eq}
\ov{r}_0(v+2\pi i)=(\phi\ot\id)(\ov{r}_0(v)).
\end{equation}
They also show that $\phi$ is a Coxeter automorphism.
The next lemma shows that in the case of trigonometric solutions coming from a solution of the AYBE
the automorphism of the Dynkin diagram is always trivial.

\begin{lem}\label{inner-lem} 
Let $r(u,v)$ be a nondegenerate unitary solution of the AYBE with the Laurent expansion 
\eqref{Laur-exp} at $u=0$. 
Assume that $\ov{r}_0(v)$ is a trigonometric solution of the CYBE with poles exactly
at $2\pi i\Z$. Then the automorphism $\phi$ in \eqref{trig-aut-pole-eq} is inner.
\end{lem}

\Pf . This follows immediately from Lemma \ref{pole-aut-lem}, since every algebra automorphism of
$A$ is inner.
\ed

Now let us recall the Belavin-Drinfeld classification of trigonometric solutions of the CYBE for $\splin_N$ 
corresponding to the trivial automorphism of $A_{N-1}$.
Let us denote by $\hg_0\sub\splin_N$ the subalgebra of traceless diagonal
matrices. 
For every Belavin-Drinfeld triple $(\Ga_1,\Ga_2,\tau)$ for $\wt{A}_{N-1}$
we have the corresponding series of solutions
\begin{equation}\label{BD-sol1}
\begin{array}{l}
\ov{r}_0(v)=t+\frac{1}{\exp(v)-1}(\pr\ot\pr)\sum_{0\le m<N, j-i\equiv m(N)}
\exp(\frac{mv}{N})e_{ij}\ot e_{ji}+\\
\sum_{0<m<N,k\ge1;j-i\equiv m(N),\tau^k(i,j)=(i',j')}
[\exp(-\frac{mv}{N})e_{ji}\ot e_{i'j'}-\exp(\frac{mv}{N})e_{i'j'}\ot e_{ji}],
\end{array}
\end{equation}
where $t\in\hg_0\ot\hg_0$ satisfies
\begin{equation}\label{BD-sol2}
t^{12}+t^{21}=(\pr\ot\pr)P^0,
\end{equation}
\begin{equation}\label{BD-sol3}
[\tau(\a)\ot \id+\id\ot\a]t=0, \ \ \a\in\Ga_1,
\end{equation}
where $P^0=\sum_i e_{ii}\ot e_{ii}$.
The result of Belavin and Drinfeld in \cite{BD} is that every nondegenerate unitary trigonometric solution
of the CYBE for $\splin_N$ that has poles exactly at $2\pi i \Z$ and the residue $(\pr\ot\pr)P$ at $0$, is conjugate to
$$\exp[v(b\ot 1)]\ov{r}_0(v)\exp[-v(b\ot 1)],$$
where $\ov{r}_0(v)$ is one of the solutions of the form \eqref{BD-sol1} and $b\in\splin_N$
is an infinitesimal symmetry of $\ov{r}_0$.

It is easy to see that the solution of the CYBE for $\splin_N$ obtained from the associative $r$-matrix
in Theorem \ref{main-thm}(i) for $S=[1,N]$ and $C_0(i)=i+1$ is given by the above formula with
\begin{equation}\label{diag-cycle}
t=\frac{1}{2}(\pr\ot\pr)P^0+s_C,
\end{equation}
where
$$s_C=\sum_{0<k<N,i}(\frac{1}{2}-\frac{k}{N})e_{ii}\ot e_{C^k(i),C^k(i)}\in\hg_0\wedge\hg_0.$$

The proof of the next result is almost identical to that of  Lemmas 4.19 and 4.20 in \cite{Sch}. 
Let us denote by $e_i:\hg\to\C$ the functional on diagonal matrices given by
$e_i(e_{jj})=\delta_{ij}$.

\begin{lem}\label{cycle-lem} 
Let $r(u,v)$ be a nondegenerate 
unitary solution of the AYBE with the Laurent expansion \eqref{Laur-exp},
such that $\ov{r}_0(v)$ is given by \eqref{BD-sol1}.
Then there exists a unique transitive cyclic permutation $C$ of $[1,N]$ such that
\eqref{diag-cycle} holds. Furthermore, for any $(i,i+1)\in\Ga_1$ with
$\tau(i,i+1)=(i',i'+1)$ one has $C(i)=i'$ and $C(i+1)=i'+1$ (i.e., $\tau$ is induced by $C\times C$).
\end{lem}

\Pf . We will make use of the identity
\begin{equation}\label{pr-AYBE-eq}
(\pr\ot\pr\ot\pr)(AYBE[\ov{r}_0])=0
\end{equation}
that follows from Theorem \ref{QYBE-thm}.
First, considering the projection of $AYBE[\ov{r}_0]$ to $\hg\ot\hg\ot\hg$
we get
$$(\pr\ot\pr\ot\pr)(AYBE[\frac{1}{\exp(v)-1}P^0+t])=0.$$
Using the fact that $t^{12}+t^{21}\equiv P^0\mod(\C\cdot 1\ot 1)$ this can be rewritten as
$$(\pr\ot\pr\ot\pr)(AYBE[t])=0.$$ 
Therefore, we have
\begin{equation}\label{pr-diag-AYBE-eq}
[(e_i-e_1)\ot (e_j-e_1)\ot (e_k-e_1)](AYBE[t])=0
\end{equation}
for all $i,j,k$.
Set $t=\sum_{i,j}t_{ij}e_{ii}\ot e_{jj}$. Note that $t_{ij}+t_{ji}=0$ for $i\neq j$ and $t_{ii}=\frac{1}{2}$
for all $i$. Let us set $t'_{ij}=t_{ij}-t_{1j}-t_{i1}$. Then substituting $t_{ij}=t'_{ij}+t_{1j}-t_{1i}$
into $t$ and then into \eqref{pr-diag-AYBE-eq} we deduce that
\begin{equation}\label{t-eq}
t'_{ij}t'_{ik}-t'_{jk}t'_{ij}+t'_{ik}t'_{jk}=\frac{1}{4}, \ \ 1<i,j,k\le N.
\end{equation}
As shown in the proof of Lemma 4.20 in \cite{Sch},
the above equation implies that $t'_{ij}=\pm\frac{1}{2}$ for $1<i,j\le N$, $i\neq j$,
and there is a unique complete order $\prec$ on $[2,N]$ such that
$t'_{ij}=\frac{1}{2}$ iff $i\prec j$ (for $i,j\in[2,N]$, $i\neq j$).
We define the cyclic permutation $C$ of $[1,N]$ by the condition that it sends each element
of $[2,N]$ to the next element with respect to this complete order.
As in the proof of Lemma 4.20 in \cite{Sch} this easily implies that 
$t-\frac{1}{2}P^0\equiv s_C\mod(\C\cdot 1\ot 1)$.

Next, we want to check that $\tau$ is induced by $C\times C$. Assume that $\tau(i,i+1)=(j,j+1)$ and
consider the coefficient $A_{ijk}$ with $e_{i+1,i}\ot e_{j,j+1}\ot e_{kk}$ in $AYBE[\ov{r}_0]$.
Let us denote by $\lan e_{lm}\ot e_{np}, r(v)\ran$ the coefficient with $e_{lm}\ot e_{np}$ in
$r(v)$. Then we have
\begin{equation}\label{coef-AYBE-eq}
\begin{array}{l}
A_{ijk}=\lan e_{i+1,i}\ot e_{j,j+1},r(v_{12})\ran \lan e_{ii}\ot e_{kk},r(v_{13})\ran-
\lan e_{jj}\ot e_{kk},r(v_{23})\ran \lan e_{i+1,i}\ot e_{j,j+1},r(v_{12})\ran+\\
\lan e_{i+1,i}\ot e_{k,k+1},r(v_{13})\ran \lan e_{j,j+1}\ot e_{k+1,k},r(v_{23})\ran.
\end{array}
\end{equation}
Note that we cannot have $\tau^n(j+1,j)=(i+1,i)$ since this would imply that
$\Ga_1$ (resp., $\Ga_2$) is the complement to $(j,j+1)$ (resp., $(i,i+1)$), $N$ is even,
$j-i\equiv N/2(N)$, and $\tau(l,l+1)=(l+N/2,l+1+N/2)$, in which case the nilpotency condition is not satisfied. Therefore, 
$$\lan e_{i+1,i}\ot e_{j,j+1},r(v)\ran=\exp(-\frac{v}{N}),$$
$$\lan e_{j,j+1}\ot e_{i+1,i},r(v)\ran=-\exp(\frac{v}{N}).$$
Next, we claim that the third summand in the right-hand side of \eqref{coef-AYBE-eq} is zero unless
$k=i$ or $k=j$. Indeed, $\tau$ (resp., $\tau^{-1}$) cannot be defined on both $(k,k+1)$ and
$(k+1,k)$. This leaves only two possibilities with $k\neq i$ and $k\neq j$: either 
$\tau^{n_1}(i,i+1)=(k,k+1)$ and $\tau^{n_2}(k,k+1)=(j,j+1)$, or
$\tau^{n_1}(j+1,j)=(k+1,k)$ and $\tau^{n_2}(k+1,k)=(i,i+1)$ (where $n_1,n_2>0$).
The latter case is impossible since $j\neq k$. In the former case we derive that
$\tau^{n_1+n_2}(i,i+1)=(j,j+1)$ which contradicts to our assumption that 
$\tau(i,i+1)=(j,j+1)$ (since $n_1+n_2\ge 2$).
Thus, recalling that
$$\lan e_{i+1,i}\ot e_{i,i+1}, r(v)\ran=\frac{\exp(\frac{(N-1)v}{N})}{\exp(v)-1},$$
$$\lan e_{j,j+1}\ot e_{j+1,j},r(v)\ran=\frac{\exp(\frac{v}{N})}{\exp(v)-1},$$
we can rewrite \eqref{coef-AYBE-eq} as follows:
\begin{align*}
&A_{ijk}=\exp(-\frac{v_{12}}{N})[t_{ik}-t_{jk}+\frac{\de_{ik}}{\exp(v_{13})-1}-\frac{\de_{jk}}{\exp(v_{23})-1}]\\
&-\de_{ik}\exp(\frac{v_{23}}{N})\frac{\exp(\frac{(N-1)v_{13}}{N})}{\exp(v_{13})-1}+
\de_{jk}\exp(-\frac{v_{13}}{N})\frac{\exp(\frac{v_{23}}{N})}{\exp(v_{23})-1}.
\end{align*}
Hence, 
$$\exp(\frac{v_{12}}{N})A_{ijk}=t_{ik}-t_{jk}-\de_{ik}.$$
Since $\pr\ot\pr\ot\pr(AYBE[\ov{r}_0])=0$, it follows that $A_{ijk}$ does not depend on $k$.
Therefore, 
$$t_{ik}-t_{jk}-\de_{ik}=[(e_i-e_j)\ot e_k]s_C-\frac{1}{2}(e_i+e_j,e_k)$$
does not depend on $k$ (note that $(e_i,e_j)=\de_{ij}$), i.e.,
\begin{equation}\label{cycle-comp1}
[(e_i-e_j)\ot \a]s_C=\frac{1}{2}(e_i+e_j,\a)
\end{equation}
for all roots $\a\in\Ga$.
Repeating the above argument for the coefficient with $e_{j,j+1}\ot e_{i+1,i}\ot e_{kk}$ in
$AYBE[\ov{r}_0]$ we derive that
\begin{equation}\label{cycle-comp2}
[(e_{i+1}-e_{j+1})\ot \a]s_C=\frac{1}{2}(e_{i+1}+e_{j+1},\a)
\end{equation}
for all $\a\in\Ga$. As shown in the proof of Lemma 4.20 in \cite{Sch}, \eqref{cycle-comp1}
and \eqref{cycle-comp2} imply that $C(i)=j$ and $C(i+1)=j+1$.
\ed

\begin{lem}\label{rescale-lem}
Assume that $N>1$. Then
a nondegenerate unitary solution $r(u,v)$ of the AYBE with the Laurent expansion at $u=0$
of the form \eqref{Laur-exp} such that $r_0(v)\equiv \ov{r}_0(v)\mod(\C\cdot 1\ot 1)$, is uniquely
determined by $\ov{r}_0$, up to rescaling $r(u,v)\mapsto \exp(\la uv)r(u,v)$.
\end{lem}

\Pf . This follows from the proof of Theorem 6 in \cite{P-AYBE}: one only has to observe that
$\ov{r}_0(v)$ is nondegenerate by Theorem \ref{QYBE-thm}(i), so it has rank $>2$ generically.
\ed

\noindent
{\it Proof of Theorem \ref{main-thm}(ii).} 
Let $r(u,v)$ be a nondegenerate unitary solution of the AYBE with the Laurent expansion at $u=0$
of the form \eqref{Laur-exp} such that $\ov{r}_0(v)$ is trigonometric. Changing $r(u,v)$ to
$cr(cu,c'v)$ we can assume that $\ov{r}_0(v)$ has poles exactly at $2\pi i\Z$ and
$\lim_{v\to 0}v\ov{r}_0(v)=(\pr\ot\pr)P$. 
Recall that we are allowed to change $r(u,v)$ to an equivalent solution
$$\wt{r}(u,v)=\exp[u(1\ot a)+v(b\ot 1)]r(u,v)\exp[-u(a\ot 1)-v(b\ot 1)],$$
where $a$ and $b$ are infinitesimal symmetries of $r(u,v)$
(note that $a$ and $b$ always commute by Lemma \ref{sym-lem}).
This operation changes $\ov{r}_0(v)$ to an equivalent solution in the sense of Belavin-Drinfeld
\cite{BD} and also changes $r_0(v)$ to $r_0(v)-a\ot 1+1\ot a$.
Therefore, in view of Lemma \ref{inner-lem} and of \eqref{r0-for}, 
changing $r(u,v)$ to an equivalent solution we can achieve that 
$r_0(v)\equiv\ov{r}_0(v)\mod(\C\cdot 1\ot 1)$ and
$\ov{r}_0(v)$ has the form \eqref{BD-sol1}. 
Note that in this case any infinitesimal symmetry of $\ov{r}_0(v)$ is diagonal (since it has to commute
with the corresponding Coxeter automorphism $\phi$ from \eqref{trig-aut-pole-eq}).
It remains to use Lemmas \ref{cycle-lem} and \ref{rescale-lem}.
\ed

\noindent
{\it Proof of Theorem \ref{AYBE-v}.} 
Let $r(v)$ be a nondegenerate unitary solution of the AYBE,
not depending on $u$. Then one can easily check that
$$r(u,v)=\frac{1\ot 1}{u}+r(v)$$
is also a nondegenerate unitary solution of the AYBE.
By Lemma \ref{pole-lem}, 
$r(u,v)$ (and hence, $r(v)$) has a simple pole at $v=0$ with the residue $cP$, where $c\neq 0$.
Now applying Lemma \ref{main-AYBE-lem2} we obtain  
$$s(u,v)=-\frac{1\ot 1}{u^2}+g(v)\cdot 1\ot 1,$$
where $g(v)=-\frac{c}{N}(\tr\ot\tr)(\frac{dr(v)}{dv})$.
Hence, by Theorem \ref{QYBE-thm}, $r(u,v)$ is a solution of the QYBE
and $\ov{r}(v)$ is a nondegenerate solution of the CYBE.
It is easy to see that $\ov{r}(v)$ cannot be equivalent to an elliptic or a trigonometric solution.
Indeed, if this were the case then $r(u,v)$ would have a pole of the form $u=u_0$ with $u_0\neq 0$
(in the elliptic case this follows from the explicit formulas for elliptic solutions in \cite{P-AYBE}, sec.2;
in the trigonometric case this follows from Theorem \ref{main-thm}(ii)). 
Now Proposition \ref{r0-prop} gives the required decomposition of $r(v)$. Therefore, 
$g(v)=-c^2/v^2$, which shows that $R(u,v)=(1/u+c/v)^{-1}r(u,v)$ satisfies 
unitarity condition \eqref{QYBE-un-eq}.
\ed

\begin{rem} The function of the form $r(v)=\frac{P}{v}+r$,
where $r\in A\ot A$ does not depend on $v$, is a unitary solution of the AYBE iff 
$r$ is a skew-symmetric constant solution of the AYBE for $A=\Mat(N,\C)$. Some information
about such solutions can be found in \cite{A2}, sec.2 (including the classification for $N=2$,
see Ex. 2.8 of {\it loc. cit.}).
\end{rem}

\end{document}